\documentclass[11pt]{amsart}
\usepackage{amssymb, amscd}
\usepackage{latexsym} 
\numberwithin{equation}{section}

\newcommand{\VVl}{\VV_{\lambda}}

\newcommand{\VV}{\mathbb V}

\newcommand{\PP}{\mathbb P}
\newcommand{\FF}{\mathbb F}
\newcommand{\FFb}{\overline{\mathbb F}}  
\newcommand{\EE}{\mathbb E}

 \newcommand{\CC}{\mathbb C}

\newcommand{\QQ}{\mathbb Q}
\newcommand{\ZZ}{\mathbb Z}
\newcommand{\RR}{\mathbb R}
\newcommand{\NN}{\mathbb N}

\newcommand{\A}{\mathcal A}
\newcommand{\M}{\mathcal M}
\newcommand{\Hh}{\mathcal H}
\newcommand{\s}{\mathbb S}
\newcommand{\Gal}[1]{{\rm Gal}(\overline{{#1}}/{#1})}
\newcommand{\Galp}[2]{{\rm Gal}({\overline{{#1}}}_{{#2}}/{{#1}}_{{#2}})}
\newcommand{\GalC}{\mathsf{Gal}_{\QQ}}
\newcommand{\MHSC}{\mathsf{MHS}}
\newcommand{\QQl}{{\mathbb Q}_{\ell}}
\newcommand{\QQlb}{{\overline{\mathbb Q}}_{\ell}}
\newcommand{\Trace}{{\rm Tr}}
\newcommand{\Aut}{{\rm Aut}}
\newcommand{\Sym}{{\rm Sym}}
\newcommand{\Sp}{{\rm{Sp}}}
\newcommand{\GSp}{{\rm{GSp}}}
\newcommand{\GL}{{\rm{GL}}}
\newcommand{\PGL}{{\rm{PGL}}}
\newcommand{\Eis}{{\rm{Eis}}}
\newcommand{\qq}{\mathbf q}
\newcommand{\gen}{\rm gen}
\newcommand{\Gtwo}{\mathrm{G}_2}

\newtheorem{theorem}{Theorem}[section]

\newtheorem{conjecture}[theorem]{Conjecture}
\newtheorem{definition-lemma}[theorem]{Definition-Lemma}

\theoremstyle{definition}
\newtheorem{definition}[theorem]{Definition}
\newtheorem{example}[theorem]{Example}

\newtheorem{notation}[theorem]{Notation}

\theoremstyle{remark}
\newtheorem{remark}[theorem]{Remark}

\begin{document}

\title[Siegel modular forms of degree three]{Siegel modular forms of degree three \\ and the cohomology of local systems}
\author{Jonas Bergstr\"om}
\address{Matematiska institutionen, Stockholms Universitet,
\newline 106 91 Stockholm, Sweden.}
\email{jonasb@math.su.se}
\author{Carel Faber}
\address{Institutionen f\"or matematik, KTH Royal Institute of Technology,
\newline 100 44 Stockholm, Sweden.}
\email{faber@math.kth.se}
\author{Gerard van der Geer}
\address{Korteweg-de Vries Instituut, Universiteit van
Amsterdam, 
\newline Postbus 94248, 1090 GE Amsterdam, The Netherlands.}
\email{geer@science.uva.nl}
\dedicatory{To the memory of Torsten Ekedahl}
\subjclass[2010]{11F46, 11G18, 14G35, 14J15, 14K10, 14H10}
\keywords{Siegel modular forms, moduli of abelian varieties, symplectic local systems, Euler characteristic, Lefschetz trace formula, Hecke operators, moduli of curves}

\begin{abstract} We give an explicit conjectural formula for the motivic Euler characteristic of an arbitrary symplectic local system on the moduli space $\A_3$ of principally polarized abelian threefolds. The main term of the formula is a conjectural motive of Siegel modular forms of a certain type; the remaining terms admit a surprisingly simple description in terms of the motivic Euler characteristics for lower genera. The conjecture is based on extensive counts of curves of genus three and abelian threefolds over finite fields. It provides a lot of new information about vector-valued Siegel modular forms of degree three, such as dimension formulas and traces of Hecke operators. We also use it to predict several lifts from genus 1 to genus 3, as well as lifts from $\Gtwo$ and new congruences of Harder type.
\end{abstract}

\maketitle

%%%%%%%%%%%%%%%%%%%%%%
\section{Introduction} \label{sec-intro}
It has been known since the 1950's that there is a close relationship between modular forms and the cohomology of local systems on moduli spaces of elliptic curves. The theorem of Eichler and Shimura expresses this relationship by exhibiting the vector space $S_{k+2}$ of cusp forms of even weight $k+2$ on ${\rm SL}(2,\ZZ)$ as the $(k+1,0)$-part of the Hodge decomposition of the first cohomology group of the local system $\VV_k$ on the moduli space  ${\A}_1$ of elliptic curves. Here, $\VV_k$ is the $k$-th symmetric power of the standard local system $\VV:=R^1\pi_* {\QQ}$ of rank $2$, where $\pi: {\mathcal E} \to {\A}_1$ is the universal elliptic curve.

Deligne strengthened the Eichler-Shimura theorem by proving that, for a prime $p$, the  trace of the  Hecke operator $T(p)$ on $S_{k+2}$ can be expressed in terms of the trace of the Frobenius map on the \'etale cohomology of the $\ell$-adic counterpart of $\VV_k$ on the moduli ${\A}_1 \otimes {\overline{\FF}}_p$ of elliptic curves in characteristic $p$. Using the Euler characteristic $e_c({\A}_1,\VV_k)= \sum (-1)^i [H^i_c({\A}_1,\VV_k)]$, we can formulate this concisely (in a suitable Grothendieck group) as 
$$ e_c({\A}_1,\VV_k)= -S[k+2]-1, $$
where $S[k+2]$ is the motive associated (by Scholl) to $S_{k+2}$ and where the term $-1$ is the contribution coming from the Eisenstein series.

The moduli space ${\A}_g$ of principally polarized abelian varieties carries the local system $\VV:=R^1\pi_* {\QQ}$ of rank $2g$, where $\pi: {\mathcal X}_g \to {\A}_g$ is the universal principally polarized abelian variety. The local system $\VV$ comes with a symplectic form, compatible with the Weil pairing on the universal abelian variety. For every irreducible representation of the symplectic group ${\rm GSp}_{2g}$, described by its highest weight $\lambda$, one has a local system $\VVl$, in such a way that $\VV$ corresponds to the dual of the standard representation of ${\rm GSp}_{2g}$. Faltings and Chai have extended the work of Eichler and Shimura to a relationship between the space $S_{n(\lambda)}$ of Siegel modular cusp forms of degree $g$ and weight $n(\lambda)$, and the cohomology of the local system $\VVl$ on $\A_g$. Note that the modular forms that appear are in general vector-valued and that the scalar-valued ones only occur for the (most) singular highest weights. The relationship leads us to believe that there should be a motivic equality of the form 
$$e_c({\A}_g,\VVl)=\sum (-1)^i [H^i_c({\A}_g, \VVl)]=(-1)^{g(g+1)/2} \, S[n(\lambda)]+e_{g,\rm extra}(\lambda),$$
generalizing the one above for genus~$1$. The (conjectural) element $S[n(\lambda)]$ of the Grothendieck group of motives should be associated to the space $S_{n(\lambda)}$ in a manner similar to the case $g=1$. As an element of the Grothendieck group of $\ell$-adic Galois representations, this means that the trace of a Frobenius element $F_p$ on $S[n(\lambda)]$ should be equal to the trace of the Hecke operator $T(p)$ on $S_{n(\lambda)}$. Our ambition is to make the equality above explicit. To do this, we need a method to find cohomological information. 

The moduli space ${\A}_g$ is defined over the integers. According to an idea of Weil, one can obtain information on the cohomology of a variety defined over the integers by counting its numbers of points over finite fields. Recall that a point of ${\A}_g$ over ${\FF}_q$ corresponds to a collection of isomorphism classes of principally polarized abelian varieties over ${\FF}_q$ that form one isomorphism class over ${\overline{\FF}}_q$; such an isomorphism class $[A]$ over ${\FF}_q$ is counted with a factor $1/|\Aut_{{\FF}_q}(A)|$. By counting abelian varieties over finite fields, one thus gets cohomological information about ${\A}_g$, and with an explicit formula as above, one would also get information about Siegel modular forms.

In fact, if one has a list of all isomorphism classes of principally polarized abelian varieties of dimension $g$ over a fixed ground field ${\FF}_q$ together with the orders of their automorphism groups and the characteristic polynomials of Frobenius acting on their first cohomology, then one can determine the trace of Frobenius (for that prime power $q$) on the Euler characteristic of the $\ell$-adic version of $\VVl$ for {\sl all} local systems $\VVl$ on ${\A}_g$. 

For $g=2$, the moduli space ${\A}_2$ coincides with the moduli space ${\M}_2^{ct}$ of curves of genus~$2$ of compact type and this makes the above counting feasible. Some years ago, the last two authors carried out the counting for finite fields of small cardinality. Subsequently, they tried to interpret the obtained traces of Frobenius on the Euler characteristic of $\VVl$ on ${\A}_2 \otimes {\FFb}_q$ in a motivic way and arrived at a precise conjecture for $e_{2,\rm extra}(\lambda)$. That is, a precise conjecture on the relation between the trace of the Hecke operator $T(p)$ on a space of Siegel cusp forms of degree $2$ and the trace of Frobenius on the $\ell$-adic \'etale cohomology of the corresponding local system on ${\A}_2\otimes {\FFb}_p$, see \cite{FvdG}. This is more difficult than in the case of genus~$1$, due to the more complicated contribution from the boundary of the moduli space and to the presence of endoscopy. The conjecture for genus~$2$ has been proved for regular local systems, but is still open for the non-regular case, see Section~\ref{sec-g2}. The conjecture and our calculations allow us to compute the trace of the Hecke operator $T(p)$ on all spaces of Siegel cusp forms of degree~$2$, for all primes $p\leq 37$. In \cite{BFG}, all three authors extended this work to degree~$2$ modular forms of level~$2$. Inspired by the results in genus~$2$, Harder formulated in \cite{Ha} a conjecture on a congruence between genus~$2$ and genus~$1$ modular forms determined by critical values of $L$-functions.

Our aim in this article is to generalize the work above to genus~$3$, that is, to give an explicit (conjectural) formula for the term $e_{3,\rm extra}(\lambda)$ in terms of Euler characteristics of local systems for genus~$1$ and~$2$. 

The conjecture is based on counts of points over finite fields, using the close relationship between ${\A}_3$ and the moduli space ${\M}_3$ of curves of genus~$3$, and a formula for the rank $1$ part of the Eisenstein cohomology. Both in genus~$2$ and genus~$3$, the motivic interpretation of the traces of Frobenius is made easier by the experimental fact that there are no Siegel cusp forms of low weight in level~$1$. Together, the conjecture and the calculations open a window on Siegel modular forms of degree $3$. We sketch some of the (heuristic) results.

The numerical Euler characteristic 
$$E_c({\A}_3,\VVl):=\sum (-1)^i \dim H^i_c({\A}_3, \VVl)$$ 
has been calculated for the local systems $\VVl$ and is known for the correction term $e_{3,\rm extra}(\lambda)$. This allows us to predict the dimension of the space of Siegel modular cusp forms of any given weight. So far, these dimensions are only known for scalar-valued Siegel modular forms by work of Tsuyumine, and our results agree with this. For $g=2$, the dimensions of most spaces of Siegel cusp forms were known earlier, so the dimension predictions are a new feature in genus~$3$.

Assuming the conjecture and using our counts, we can calculate the trace of $T(p)$ on any space of Siegel modular forms $S_{n(\lambda)}$ of degree $3$ for $p\leq 17$. Moreover, if $\dim S_{n(\lambda)}=1$, we can compute the local spinor $L$-factor at $p=2$. 

We make a precise conjecture on lifts from genus~$1$ to genus~$3$. In particular, for every triple $(f,g,h)$ of elliptic cusp forms that are 
Hecke eigenforms of weights $b+3$, $a+c+5$, and $a-c+3$,
where $a\ge b\ge c\ge 0$,
there should be a Siegel modular cusp form $F$ of weight $(a-b,b-c,c+4)$ that is an eigenform for the Hecke algebra with spinor $L$-function $L(f\otimes g,s)L(f\otimes h, s-c-1)$.

We find strong evidence for the existence of Siegel modular forms of degree~$3$ (and level $1$) that are lifts from $\Gtwo$, as predicted by Gross and Savin. Finally, we are able to formulate conjectures of Harder type on congruences between Siegel modular forms of degree $1$ and degree $3$. All these heuristic findings provide strong consistency checks for our main conjecture.

We hope that our results may help to make (vector-valued) Siegel modular forms of degree~$3$ as concrete as elliptic modular forms.

Our methods apply to some extent also to local systems on the moduli space ${\M}_3$ and Teichm\"uller modular forms; for more on this we refer to the forthcoming article \cite{BFG2}.

After reviewing the case of genus~$1$, we introduce Siegel modular forms and give a short summary of the results of Faltings and Chai that we shall use. We then discuss the hypothetical motive attached to the space of cusp forms of a given weight. In Section \ref{sec-g2}, we review and reformulate our results for genus~$2$ in a form that is close to our generalization for genus~$3$. We present our main conjectures in Section \ref{sec-g3}. The method of counting is explained in Section \ref{sec-points}. We then discuss the evidence for Siegel modular forms of type $\Gtwo$. In the final section, we present our Harder-type conjectures on congruences between Siegel modular forms of degree $3$ and degree $2$ and $1$.

This work may be viewed as belonging to the Langlands program in that we study the connection between automorphic forms on $\GSp_{2g}$ and Galois representations appearing in the cohomology of local systems on $\A_g$. Instead of using the Lefschetz trace formula, as we do here, one can find cohomological information of this kind via the Arthur-Selberg trace formula. The study of the latter trace formula is a vast field of research; the reference \cite{Arthur} could serve as an introduction. 
The group ${\rm GSp}_{2g}$ has been investigated by several people, see for instance the references \cite{Kottwitz, Laumon,Morel,Weissauer-book}. 
Clearly, the Arthur-Selberg trace formula has been used successfully in these and other references; to our knowledge, however, this method does not produce formulas as explicit as the ones in this article.
We hope that researchers in automorphic forms will be interested in our results, even though we use a rather different language.
See \cite{C-R} for some very interesting recent developments
closely related to our work.

%%%%%%%%%%%%%%%%%%%%%%%%%%%
\section*{Acknowledgements}
The authors thank Pierre Deligne, Neil Dummigan, Benedict Gross, G\"un\-ter Harder, Anton Mellit, and Don Zagier for their contributions, and the Max Planck Institute for Mathematics in Bonn for hospitality and excellent working conditions. We are very grateful to Maarten Hoeve for assistance with the computer programming and we thank Hidenori Katsurada for his remarks about congruences. 
Finally, we thank the referee.

The second author was supported by
the G\"oran Gustafsson Foundation for Research in Natural Sciences
and Medicine (KVA) and grant 622-2003-1123 from the Swedish Research Council.
The third author's visit to KTH in October 2010 was supported by
the G\"oran Gustafsson Foundation for Research in Natural Sciences
and Medicine (UU/KTH).

%%%%%%%%%%%%%%%%%%%
\section{Genus one} \label{sec-g1} 
We start by reviewing the theorem in genus~$1$ by Eichler, Shimura and Deligne that we wish to generalize to higher genera. This theorem describes the cohomology of certain local systems on the moduli space of elliptic curves in terms of elliptic modular forms. 

Let $\pi\colon{\mathcal X}_1 \to \A_1$ be the universal object over the moduli space of elliptic curves. These spaces are smooth Deligne-Mumford stacks over $\ZZ$. 

First, we consider the analytic picture. Define a local system $\VV:=R^1 \pi_*\CC$ on $\A_1 \otimes \CC$. For any $a \geq 0$ we put $\VV_a:=\Sym^a (\VV)$. We wish to understand the compactly supported Betti cohomology groups $H_c^i(\A_1 \otimes \CC,\VV_a)$ with their mixed Hodge structures. Define the inner cohomology group $H_!^i(\A_1 \otimes \CC,\VV_a)$ as the image of the natural map
$$H_c^i(\A_1 \otimes \CC,\VV_a) \to H^i(\A_1 \otimes \CC,\VV_a),$$
and the Eisenstein cohomology group $H_{\Eis}^i(\A_1 \otimes \CC,\VV_a)$ as the kernel of the same map. Since the element $-{\rm Id}$ of ${\rm SL}(2,\ZZ)$ acts as $-1$ on $\VV$, all these cohomology groups vanish if $a$ is odd. We therefore assume from now on that $a$ is even. 

The group ${\rm SL}(2,\ZZ)$ acts on the complex upper half space $\mathfrak{H}_1$ by $z \mapsto \alpha(z):=(az+b)(cz+d)^{-1}$ for any $z \in \mathfrak{H}_1$ and $\alpha=(a,b\,;c,d) \in {\rm SL}(2,\ZZ)$. We note that $\A_1\otimes \CC \cong {\rm SL}(2,\ZZ) \backslash \mathfrak{H}_1$ as analytic spaces. We define an elliptic modular form of weight $k$ as a holomorphic map $f: \mathfrak{H}_1 \rightarrow \CC$, such that $f$ is also holomorphic at the point in infinity and such that $f(\alpha(z))=(cz+d)^k f(z)$ for all $z \in \mathfrak{H}_1$ and $\alpha \in {\rm SL}(2,\ZZ)$. We call an elliptic modular form a \emph{cusp} form if it vanishes at the point in infinity. Let $S_k$ be the vector space of elliptic cusp forms of weight $k$ and put $s_k:=\dim_{\CC} S_k$. The \emph{Hecke algebra} acts on $S_k$. It is generated by operators $T(p)$ and $T_1(p^2)$ for each prime $p$. These operators are simultaneously diagonalizable and the eigenvectors will be called \emph{eigenforms}. The eigenvalue for $T(p)$ of an eigenform $f$ will be denoted by $\lambda_p(f)$. 

\begin{example} The best known cusp form is probably $\Delta$. It has weight $12$ and can be defined by $$\Delta(q):=q \prod_{n=1}^{\infty}(1-q^n)^{24}=\sum_{n=1}^{\infty} \tau(n) q^n,$$ 
where $q:=e^{2 \pi i z}$ and $\tau$ is the Ramanujan tau function. Since $s_{12}=1$, $\Delta$ is an eigenform and $\lambda_p(\Delta)=\tau(p)$ for every prime $p$.
\end{example}

Similarly to the above, we define a bundle $\EE:=\pi_* \Omega_{{\mathcal X}_1 / \A_1}$ on $\A_1 \otimes \CC$ and for any $a \geq 0$ we put $\EE_a:=\Sym^a (\EE)$. The moduli space $\A_1$ can be compactified by adding the point at infinity, giving $\A'_1:=\A_1 \cup \{\infty\}$. The line bundles $\EE_a$ can be extended to $\A'_1$. We have the identification 
$$H^0\bigl(\A'_1 \otimes \CC,\EE_{k}(-\infty)\bigr)\cong S_k.$$ 
The Eichler-Shimura theorem (see \cite[Th\'eor\`eme 2.10]{Deligne}) then tells us that
$$S_{a+2} \oplus \overline{S}_{a+2} \cong H_!^1(\A_1 \otimes \CC, \VV_a),$$
where $S_{a+2}$ has Hodge type $(a+1,0)$. If $a>0$, then $H_{\Eis}^1(\A_1 \otimes \CC,\VV_a) \cong \CC$ with Hodge type $(0,0)$, whereas $H_{c}^0(\A_1 \otimes \CC,\VV_a)$ and $H_{c}^2(\A_1 \otimes \CC,\VV_a)$ vanish. 

In order to connect the arithmetic properties of elliptic modular forms to these cohomology groups, we turn to $\ell$-adic \'etale cohomology and we redo our definitions of the local systems by putting $\VV:=R^1 \pi_* \QQ_{\ell}$ and $\VV_a:=\Sym^a (\VV)$, which are $\ell$-adic local systems over $\A_1$. For each prime $p$, we define a geometric Hecke operator, also denoted by $T(p)$, by using the correspondence in characteristic $p$ coming from the two natural projections to $\A_1$ from the moduli space of cyclic $p$-isogenies between elliptic curves.
The operators $T(p)$ act on $H_c^i(\A_1 \otimes \FF_p,\VV_a)$
and $H^i(\A_1 \otimes \FF_p,\VV_a)$. The geometric action of $T(p)$ on $H_!^1(\A_1 \otimes \FF_p,\VV_a)$ is then equal to the action of $T(p)$ on $S_{a+2} \oplus \overline{S}_{a+2}$, see \cite[Prop.~3.19]{Deligne}. 

On the other hand, the cohomology groups $H^i_c(\A_1 \otimes \overline{K},\VV_a)$ come with an action of $\Gal{K}$. For any prime $p$, there is a natural 
surjection $\Galp{\QQ}{p} \to \Galp{\FF}{p}$, and if $p \neq \ell$, there is an isomorphism
$$H^i_c(\A_1 \otimes \FFb_p,\VV_a) \to H^i_c(\A_1 \otimes \overline{\QQ}_p,\VV_a)$$
of $\Galp{\QQ}{p}$-representations. This isomorphism also holds for Eisenstein and inner cohomology. We define the (geometric) Frobenius map $F_q \in \Galp{\FF}{q}$ to be the inverse of $x \mapsto x^q$. The two actions are connected in the following way:
$$\Trace\bigl(F_p,H^1_!(\A_1 \otimes \FFb_p,\VV_a)\bigr)=\Trace(T(p),S_{a+2}),$$
see \cite[Prop.~4.8]{Deligne}. We also note that for $a>0$, 
$$\Trace\bigl(F_p,H^1_{\Eis}(\A_1 \otimes \FFb_p,\VV_a)\bigr)=1$$
for all primes $p \neq \ell$. 

We can choose an injection $\Galp{\QQ}{p} \to \Gal{\QQ}$ and talk about Frobenius elements in $\Gal{\QQ}$ as elements of $\Galp{\QQ}{p}$ that are mapped to $F_p \in \Galp{\FF}{p}$. We will, for each prime $p$, choose such a Frobenius element and it will by abuse of notation also be denoted by $F_p$. 

\begin{remark} \label{rmk-traces} The traces of $F_p$ for all (unramified) primes $p$ determine an $\ell$-adic $\Gal{\QQ}$-representation up to semi-simplification, see \cite[Proposition 2.6]{Fermat}.
\end{remark}

For any $p \neq \ell$, the two $\Galp{\QQ}{p}$-representations $H^i_c(\A_1 \otimes \overline{\QQ}_p,\VV_a)$ and $H^i_c(\A_1 \otimes \overline{\QQ},\VV_a)$ are isomorphic, and the same holds for inner cohomology. It follows that 
$$\Trace\bigl(F_p,H^1_!(\A_1 \otimes \overline{\QQ},\VV_a)\bigr)=\Trace(T(p),S_{a+2}),$$
and this then determines the $\Gal{\QQ}$-representation $H^1_!(\A_1 \otimes \overline{\QQ},\VV_a)$ up to semi-simplification. 

For $a \geq 2$, there is a construction by Scholl \cite[Theorem 1.2.4]{Scholl} of a corresponding Chow motive $S[a+2]$. This motive is defined over $\QQ$, it has rank $2 \, s_{a+2}$, its Betti realization has a pure Hodge structure with types $(a+1,0)$ and $(0,a+1)$, and its $\ell$-adic realization has the property that 
$$\Trace(F_p,S[a+2])=\Trace(T(p),S_{a+2}),$$ 
for all primes $p \neq \ell$. See also \cite{C-F} for an alternative construction of $S[a+2]$. 

We can then rewrite the results above in terms of a motivic Euler characteristic:
$$e_c({\A}_1\otimes \QQ,\VV_a):=\sum_{i=0}^2 (-1)^i \, [H^i_c(\A_1 \otimes \QQ,\VV_a)].$$ 
Here follows the theorem that we wish to generalize to higher genera. 
\begin{theorem} \label{thm-g1} For every even $a >0$,
$$e_c({\A}_1\otimes \QQ,\VV_a)=-S[a+2]-1.$$ 
\end{theorem}

\begin{remark} As a result, one can compute $\Trace(T(p),S_{a+2})$ by means of $\Trace(F_p,e_c({\A}_1\otimes \FF_p,\VV_a))$, which in turn can be found by counting elliptic curves over $\FF_p$ together with their number of points over $\FF_p$. Essentially, make a list of all elliptic curves defined over $\FF_p$ up to $\FF_p$-isomorphism; determine for each $E$ in this list $|\Aut_{\FF_p}(E)|$ and $|E(\FF_p)|$. Having done this, one can calculate $\Trace(T(p),S_{a+2})$ for all $a >0$. But there are of course other (possibly more efficient) ways of computing these numbers, see for instance the tables of Stein \cite{Stein}. 
\end{remark}

For $a=0$, we have $e_c({\A}_1 \otimes \QQ,\VV_0)=L$, where $L$ is the so called \emph{Lefschetz motive}. The $\ell$-adic realization of $L$ equals the cyclotomic representation $\QQ_{\ell}(-1)$, which is one-dimensional and satisfies $\Trace(F_q,\QQl(-1))=q$ for all prime powers $q$. For bookkeeping reasons, we want the formula in Theorem~\ref{thm-g1} to hold also for $a=0$, so we define $S[2]:=-L-1$. To be consistent, we then also put $s_2:=-1$.

For further use below, we note that in \cite{Scholl} there is actually a construction of a motive $M_f$ for any eigenform $f \in S_k$. This motive will be defined over some number field and its $\ell$-adic realization has the property that $\Trace(F_p,M_f)=\Trace(T(p),f)=\lambda_p(f)$.

%%%%%%%%%%%%%%%%%%%%%%%%%%%%%%
\section{Siegel modular forms} \label{sec-smf} 
In this section, we recall the notion of Siegel modular forms of degree $g$, which are natural generalizations of the elliptic modular forms that occur in genus~$1$. General references are \cite{Andrianov}, \cite{F-C}, \cite{Freitag}, and \cite{vdG1}.

The group $\Sp(2g,\ZZ)$ acts on the Siegel upper half space 
$$\mathfrak{H}_g:=\{z \in {\rm Mat}(g\times g): z^{t}=z, 
{\rm Im}(z) > 0\}$$ 
by $\tau \mapsto \alpha(\tau):=(a\tau+b)(c\tau+d)^{-1}$ for any $\tau \in \mathfrak{H}_g$ and $\alpha=(a,b;c,d) \in \Sp(2g,\ZZ)$. Let $W$ be a finite-dimensional complex vector space and $\rho: \GL(g,\CC) \rightarrow \GL(W)$ an irreducible representation. A \emph{Siegel modular form} of degree $g$ and weight $\rho$ is a holomorphic map $f: \mathfrak{H}_g \rightarrow W$ such that $f(\alpha(\tau))=\rho(c\tau+d)f(\tau)$ for all $\tau \in \mathfrak{H}_g$ and $\alpha \in \Sp(2g,\ZZ)$. When $g=1$, we also require $f$ to be holomorphic at infinity. Let $U$ be the standard representation of $\GL(g,\CC)$. For each $g$-tuple $(n_1,\ldots,n_g) \in \NN^g$, we define $U_{(n_1,\ldots,n_g)}$ to be the irreducible representation of $\GL(g,\CC)$ of highest weight in 
$$\Sym^{n_1}(\wedge^1 U) \otimes \Sym^{n_2}(\wedge^2 U) \otimes \ldots \otimes \Sym^{n_{g-1}}(\wedge^{g-1} U) \otimes (\wedge^{g} U)^{n_{g}}.$$ 
It can be cut out by Schur functors. We have for instance $U \otimes \wedge^2 U \cong U_{1,1,0} \oplus U_{0,0,1}$ for $g=3$. We let the cotangent bundle $E:=\pi_* \Omega^1_{\mathcal{X}_g/\A_g}$ (the Hodge bundle) correspond to the standard representation of $\GL(g,\CC)$ and using the above construction we get vector bundles $E_{(n_1,\ldots,n_g)}$. For $g > 1$, we can then identify the vector space of Siegel modular forms of weight $(n_1,\ldots,n_g)$ with $H^0(\A_g \otimes \CC,E_{(n_1,\ldots,n_g)})$. If we take any Faltings-Chai toroidal compactification $\A'_g$ of $\A_g$ and if we let $D:=\A'_g \setminus \A_g$ be the divisor at infinity, then we can define the vector space of Siegel modular \emph{cusp} forms, $S_{(n_1,\ldots,n_g)}$, to be $H^0(\A'_g \otimes \CC,E_{(n_1,\ldots,n_g)}(-D))$. In other words, a Siegel modular form is a cusp form if it vanishes along the divisor at infinity. We will call a Siegel modular cusp form \emph{classical} if it is scalar-valued, i.e., if $n_1=n_2=\ldots=n_{g-1}=0$. Let us also put $s_{(n_1,\ldots,n_g)}:=\dim_{\CC} S_{(n_1,\ldots,n_g)}$.

The \emph{Hecke algebra}, whose elements are called \emph{Hecke operators}, acts on the space $S_{(n_1,\ldots,n_g)}$. It is a tensor product over all primes $p$ of \emph{local} Hecke algebras that are generated by elements of the form $T(p)$ and $T_i(p^2)$ for $i=1,\ldots,g$, see \cite{Andrianov} or \cite{Freitag}. These operators are simultaneously diagonalizable and we call the eigenvectors eigenforms with corresponding eigenvalues $\lambda_p$ and $\lambda_{i,p^2}$. To a Hecke eigenform we can then associate a homomorphism from the Hecke algebra to $\CC$.

The Satake isomorphism identifies the local Hecke algebra at any prime~$p$ with the representation ring of $\mathrm{GSpin}_{2g+1}(\CC)$, the dual group of $\GSp_{2g}$. Thus, a Hecke eigenform $f$ determines, for each prime $p$, a conjugacy class $s_p(f)$ in $\mathrm{GSpin}_{2g+1}(\CC)$. If we fix a representation $r$ of $\mathrm{GSpin}_{2g+1}(\CC)$, we can form an $L$-function by letting the local factor at $p$ equal $Q_p(p^{-s},f)^{-1}$, where $Q_p(X,f)$ is the characteristic polynomial $\det(1-r(s_p(f))X)$ of $r(s_p(f))$, see \cite[p.~50]{Borel}. In this article, we will let $r$ be the spin representation, from which we get the so called \emph{spinor} $L$-function $L(f,s)$.

The local Hecke algebra of $\GSp_{2g}$ can also be identified with the elements of the local Hecke algebra of the diagonal torus in $\GSp_{2g}$ that are fixed by its Weyl group. Using this, we can associate to a Hecke eigenform $f$ the $(g+1)$-tuple of its Satake $p$-parameters $(\alpha_0(f),\alpha_1(f),\ldots,\alpha_g(f)) \in \CC^{g+1}$. The local factor $L_p(s,f)$ of the spinor $L$-function of $f$ then equals $Q_p(p^{-s},f)^{-1}$, where 
$$Q_p(X,f):=\prod_{r=0}^g\, \prod_{1\leq i_1 < i_2 < \ldots <i_r\leq g}\bigl(1-\alpha_0(f) \alpha_{i_1}(f) \cdots \alpha_{i_r}(f)X \bigr).$$
Note that the polynomial $Q_p(X,f)$ has degree $2^g$ and the coefficient of $X$ is equal to $-\lambda_p(f)$.

\begin{example} \label{exa-miyawaki} In \cite[p.~314]{Miyawaki}, Miyawaki constructed a non-zero classical cusp form $F_{12} \in S_{0,0,12}$. Define 
$$E_8:=\{(x_1,\ldots,x_8)\in \RR^8: 2x_i\in \ZZ, x_i-x_j\in \ZZ, x_1+\ldots+x_8\in 2\ZZ \},$$ 
which is the unique even unimodular lattice of rank $8$. Let $I_3$ be the $3 \times 3$ identity matrix and define $Q$ to be the $3 \times 8$ matrix $(I_3,i \cdot I_3,0,0)$. If $\langle \cdot,\cdot \rangle$ denotes the usual inner product on $\RR^8$, then 
$$F_{12}(\tau):=\sum_{v_1,v_2,v_3 \in E_8} {\rm{Re}}\Bigl({\rm det}\bigl(Q(v_1,v_2,v_3)\bigr)^8\Bigr){\rm exp}\Bigl(\pi i \, \Trace\bigl((\langle v_i,v_j \rangle)\tau\bigr)\Bigr).$$
Since $s_{0,0,12}=1$, $F_{12}$ is an eigenform and Miyawaki computed the eigenvalues for $T(2)$, $T_1(4)$, $T_2(4)$, and $T_3(4)$, where the first equals $-2^6 \cdot 747$. Based on these computations, he conjectured for all primes $p$ the equality
$$\lambda_p(F_{12})=\lambda_p(\Delta) \bigl(\lambda_p(f)+p^9+p^{10} \bigr),$$ 
where $f$ is an eigenform in $S_{20}$. This was proved by Ikeda, see further in Section~\ref{sec-g3-nonexh}. 
\end{example}

%%%%%%%%%%%%%%%%%%%%%%%%%%%%%%%%%%%%%
\section{Cohomology of local systems} \label{sec-localsystems}
In this section we introduce the Euler characteristic that we would like to compute and review results of Faltings and Chai on the cohomology in question.

Let $\M_g$ be the moduli space of smooth curves of genus~$g$ and $\A_g$ the moduli space of principally polarized abelian varieties of dimension $g$. These are smooth Deligne-Mumford stacks defined over ${\rm Spec}(\ZZ)$. 

Using the universal abelian variety $\pi\colon{\mathcal X}_g \to \A_g$, we define a local system on $\A_g$. It comes in a Betti version ${}_0\VV:= R^1\pi_* \QQ$ on $\A_g \otimes \QQ$ and an $\ell$-adic version ${}_{\ell}\VV:=R^1\pi_* \QQl$ on $\A_g \otimes \ZZ[1/\ell]$. We will often denote both of them simply by $\VV$. 

For any $[A] \in \A_g$, the stalk $(\VV)_A$ is isomorphic to $H^1(A)$ (with coefficients in $\QQ$ or $\QQl$ depending on the cohomology theory). Using the polarization and Poincar\'e duality, we get a symplectic pairing ${}_0 \VV \times {}_0 \VV \rightarrow \QQ(-1)$, where  $\QQ(-1)$ is a Tate twist (and similarly, ${}_{\ell} \VV \times {}_{\ell} \VV \rightarrow \QQl(-1)$, where $\QQl(-1)$ corresponds to the inverse of the cyclotomic character). 

Let $V$ be the contragredient of the standard representation of $\GSp_{2g}$, which is isomorphic to the tensor product of the standard representation of $\GSp_{2g}$ and the inverse of the multiplier representation $\eta$, see \cite[p.~224]{F-C}. We will consider partitions $\lambda$ of length at most $g$ and they will be written in the form $\lambda=(\lambda_1 \geq \lambda_2 \geq \ldots \geq \lambda_g \geq 0)$. For each such partition $\lambda$, we define the representation $V_{\lambda}$ of $\GSp_{2g}$ to be the irreducible representation of highest weight in 
$$\Sym^{\lambda_1-\lambda_2}(\wedge^1 V) 
\otimes \ldots 
\otimes \Sym^{\lambda_{g-1}-\lambda_g}(\wedge^{g-1} V) 
\otimes \Sym^{\lambda_{g}}(\wedge^{g} V).$$
It can be cut out using Schur functors. 
This gives all irreducible representations of $\GSp_{2g}$ modulo tensoring with $\eta$. 
For instance, for genus~$2$ we have
$\wedge^2 V \cong V_{1,1} \oplus V_{0,0}\otimes\eta^{-1}$.

Our local system $\VV$ corresponds now to the equivariant bundle defined by the contragredient of the standard representation of $\GSp_{2g}$. By applying the construction above to $\VV$, we define a local system $\VVl$. For example, we have $\wedge^2 \VV \cong \VV_{1,1} \oplus \VV_{0,0}(-1)$ for genus~$2$.
If $\lambda_1 > \ldots > \lambda_g >0$, then this local system is called $\emph{regular}$. For a partition $\lambda$, we define $|\lambda|:=\sum_{i=1}^g \lambda_i$ to be the weight of $\lambda$. 

We are interested in the \emph{motivic} Euler characteristic 
\begin{equation}\label{eq-euler} 
e_c(\A_g,\VVl)= \sum_{i=0}^{g(g+1)} (-1)^i \, [H_c^i(\A_g,\VVl)]. 
\end{equation}
In practice, we will be considering this expression either in the Grothendieck group of mixed Hodge structures (denoted $K_0(\MHSC)$), by taking the compactly supported Betti cohomology of ${}_0\VVl$ on $\A_g \otimes \CC$, or in the Grothendieck group of $\ell$-adic $\Gal{\QQ}$-representations (denoted $K_0(\GalC)$), by taking the compactly supported $\ell$-adic \'etale cohomology of ${}_{\ell}\VVl$ on $\A_g \otimes \overline{\QQ}$. 

\begin{remark} Note that further tensoring with $\eta$ does not give new interesting local systems, since it corresponds to Tate twists, that is, $H_c^i(\A_g,\VVl(-j)) \cong H_c^i(\A_g,\VVl)(-j)$. 
\end{remark}

\begin{remark} The element $-{\rm Id}$ belongs to $\GSp_{2g}$ and it acts as $-1$ on $\VV$. This has the consequence that $H_c^*(\A_g,\VVl)=0$ if $|\lambda|$ is odd. From now on, we therefore always assume that $|\lambda|$ is even.
\end{remark} 

We also define the \emph{integer}-valued Euler characteristic:
\begin{equation*} 
E_c(\A_g,\VVl)= \sum_{i=0}^{g(g+1)} (-1)^i  \dim \bigl(H_c^i(\A_g \otimes \CC,{}_0 \VVl) \bigr).
\end{equation*}

\subsection{Results of Faltings and Chai} \label{sec-FC}
The cohomology groups $H_c^i(\A_g\otimes \CC,\VVl)$ and $H^i(\A_g\otimes \CC,\VVl)$ carry mixed Hodge structures of weight $\leq |\lambda|+i$ respectively $\geq |\lambda|+i$, see \cite[p.~233]{F-C}. Since ${}_{\ell}\VVl$ is a sheaf of pure weight $|\lambda|$, the same weight claim holds for the $\ell$-adic cohomology in the sense of Deligne, see \cite[Cor.~3.3.3, 3.3.4]{DeligneWII}. The steps in the Hodge filtration for the cohomology groups are given by the sums of the elements of any of the $2^g$ subsets of $\{\lambda_g+1,\lambda_{g-1}+2,\ldots,\lambda_1+g\}$. In genus~$3$, the explicit Hodge filtration for $\lambda=(a\ge b\ge c\ge 0)$ is
$$F^0  \supseteq F^{c+1} \supseteq  F^{b+2} \supseteq F^{t_2} \supseteq F^{t_1} \supseteq F^{a+c+4} \supseteq F^{a+b+5} \supseteq F^{a+b+c+6},$$
where $t_1 \geq t_2$ and $\{t_1,t_2\}=\{b+c+3,a+3\}$. For $H_c^{\bullet}(\A_3\otimes \CC,\VVl)$ and $a \neq b+c$, the graded pieces in the Hodge filtration can be identified with the following coherent cohomology groups:
$$
\begin{aligned}
F^0/F^{c+1} \cong &\,H^{\bullet-0}({\A'_3}\otimes \CC,E_{b-c,a-b,-a}(-D)), \\ 
F^{c+1}/F^{b+2} \cong &\,H^{\bullet-1}({\A'_3}\otimes \CC,E_{b+c+2,a-b,-a}(-D)) ,\\
F^{b+2}/F^{r_1} \cong &\,H^{\bullet-2}({\A'_3}\otimes \CC,E_{b+c+2,a-c+1,-a}(-D)), \\ 
F^{a+3}/F^{r_2} \cong &\,H^{\bullet-3}({\A'_3}\otimes \CC,E_{a+c+3,b-c,1-b}(-D)),\\
F^{b+c+3}/F^{r_3} \cong  &\,H^{\bullet-3}({\A'_3}\otimes \CC,E_{b-c,a+c+3,-a}(-D)), \\ 
F^{a+c+4}/F^{a+b+5} \cong &\,H^{\bullet-4}({\A'_3}\otimes \CC,E_{a-c+1,b+c+2,1-b}(-D)),\\ 
F^{a+b+5}/F^{a+b+c+6} \cong &\,H^{\bullet-5}({\A'_3}\otimes \CC,E_{a-b,b+c+2,2-c}(-D)), \\ 
F^{a+b+c+6} \cong &\,H^{\bullet-6}({\A'_3}\otimes \CC,E_{a-b,b-c,c+4}(-D)), \\
\end{aligned}
$$ 
where $r_1,r_2,r_3$ depend in the obvious way on the ordering of $b+c+3$ and $a+3$. If $r=a+3=b+c+3$, then the above holds, except that $F^r/F^{a+c+4} \cong  H^{\bullet-3}({\A'_3}\otimes \CC,E_{a+c+3,b-c,1-b}(-D))\oplus H^{\bullet-3}({\A'_3}\otimes \CC,E_{b-c,a+c+3,-a}(-D))$.

\begin{notation} For any partition $\lambda$, we put 
$$n(\lambda):=(\lambda_1-\lambda_2,\lambda_2-\lambda_3,\ldots,\lambda_{g-1}-\lambda_g,\lambda_g+g+1).$$ 
\end{notation}

The last step of the Hodge filtration of $H_c^{g(g+1)/2}(\A_g \otimes \CC,\VVl)$ is given by the Siegel modular cusp forms of weight $n(\lambda)$, that is, 
$$F^{|\lambda|+g(g+1)/2} \cong H^{0}({\A'_g}\otimes \CC,E_{n(\lambda)}(-D)) \cong S_{n(\lambda)},$$ see \cite[p.~237]{F-C}.  

We define the inner cohomology $H_!^i(\A_g,\VVl)$ as the image of the natural map $H^i_c(\A_g,\VVl) \to H^i(\A_g,\VVl)$. It is pure of weight $|\lambda|+i$. Define the Eisenstein cohomology $H_{\Eis}^i(\A_g,\VVl)$ as the kernel of the same map and let $e_{g,\Eis}(\lambda)$ denote the corresponding Euler characteristic. If $\lambda$ is regular, then $H_!^i(\A_g,\VVl) = 0$ for $i \neq g(g+1)/2$, see \cite{Faltings}.

%%%%%%%%%%%%%%%%%%%%%%%%%%%%%%%%%%%%%%%%%%%%
\section{The motive of Siegel modular forms} \label{sec-motives} 
The formula in genus~$1$ (for $a$ even and positive) 
\begin{equation} \label{eq-g1} 
e_c(A_1 \otimes {\QQ},\VV_a)=-1-S[a+2] 
\end{equation}
of Theorem~\ref{thm-g1} and the results of Faltings and Chai (see Section~\ref{sec-FC}) suggest (or invite) a generalization to higher $g$. Unfortunately, the generalization is not straightforward. On the one hand, the Eisenstein cohomology (the analogue of the $-1$ in Equation~\eqref{eq-g1}) is more complicated, and on the other hand, there are contributions from endoscopic groups. Furthermore, it is not known how to define the analogues of $S[a+2]$ for higher $g$. But the first expectation is that each Siegel modular form of degree $g$ and weight $n(\lambda)$ that is an eigenform of the Hecke algebra should contribute a piece of rank $2^g$ to the middle inner cohomology group. We therefore introduce $S[n(\lambda)]$, a conjectural element of the \emph{Grothendieck group} of motives defined over $\QQ$, of rank $2^g \, s_{n(\lambda)}$ and whose $\ell$-adic realization should have the following property for all primes $p \neq \ell$: 
$$\Trace(F_p,S[n(\lambda)])=\Trace(T(p),S_{n(\lambda)}).$$
This property determines $S[n(\lambda)]$ as an element of  $K_0(\GalC)$, see Remark~\ref{rmk-traces}. 

The (conjectural) Langlands correspondence connects, loosely speaking, automorphic forms of a reductive group $G$ with continuous homomorphisms from the Galois group to the dual group $\widehat G$. In our case, the spinor $L$-function of a Siegel modular form for $\GSp_{2g}$ (see Section~\ref{sec-smf}) should equal the $L$-function of a Galois representation 
$$\Gal{\QQ} \rightarrow  \widehat {\GSp}_{2g}(\QQlb) \cong \mathrm{GSpin}_{2g+1}(\QQlb) \rightarrow \mathrm{GL}(2^g,\QQlb),$$ 
where the last arrow is the spin representation. That is, if $f_1,\ldots, f_{s_{n(\lambda)}}$ is a basis of Hecke eigenforms of $S_{n(\lambda)}$, then the characteristic polynomial $C_p(X,S[n(\lambda)])=\det(1-F_pX)$ of Frobenius acting on the $\ell$-adic realization of $S[n(\lambda)]$ should equal the product of the characteristic polynomials of the corresponding Hecke eigenforms, i.e., 
$$C_p(X,S[n(\lambda)])=\prod_{i=1}^{s_{n(\lambda)}} Q_p(X,f_i).$$

As said above, the first expectation is that each Hecke eigenform contributes a piece of rank $2^g$ to the middle inner cohomology group. However, this expectation fails: some modular forms contribute a smaller piece; these will be called \emph{non-exhaustive}. So we also introduce $\widehat{S}[n(\lambda)]$, another conjectural element of the Grothendieck group of motives over $\QQ$. It should correspond to the direct sum of the \emph{actual} contributions to the middle inner cohomology group coming from the various Hecke eigenforms. We will continue to use $S[n(\lambda)]$ as well; it is surprisingly useful as a bookkeeping device.

We will say that a direct summand $\Sigma_{n(\lambda)}$ of $S_{n(\lambda)}$ as a Hecke module over $\QQ$ has the expected properties if there is a submotive $\Sigma[n(\lambda)]$ of $\widehat{S}[n(\lambda)]$ of rank $2^g \, \dim \Sigma_{n(\lambda)}$ such that if $f_1,\ldots f_{k}$ is a basis of Hecke eigenforms of $\Sigma_{n(\lambda)}$, then 
$$C_p(X,\Sigma[n(\lambda)])=\prod_{i=1}^{k} Q_p(X,f_k)$$
holds for all $p \neq \ell$; moreover, as an element of $K_0(\MHSC)$, the dimension of the piece of Hodge type $(r,|\lambda|+g(g+1)/2-r)$ should equal $\dim \Sigma_{n(\lambda)}$ times the number of subsets of the list $(\lambda_g+1,\lambda_{g-1}+2,\ldots,\lambda_1+g)$ that have sum $r$. If the whole of $S_{n(\lambda)}$ has the expected properties (which happens when there are no non-exhaustive forms), then we will call $\lambda$ \emph{normal}.

In genus~$1$, the motive $S[a+2]$ for $a \geq 2$ has been constructed. It appears in the first inner cohomology group of $\VV_{(a)}$ on $\A_1$ and $\lambda=(a)$ is normal for all $a \geq 2$. As to $a=0$, the inner cohomology of $\A_1$ with $\QQl$ coefficients vanishes and thus $\widehat{S}[2]=0$; purely for bookkeeping reasons, we earlier defined $S[2]=-L-1$ and thus $s_2=-1$.

The inner cohomology does not only consist of $\widehat{S}[n(\lambda)]$; what is left should be contributions connected to the so called \emph{endoscopic groups}. We call these contributions the \emph{endoscopic cohomology}, and we denote its Euler characteristic by $e_{g,{\rm endo}}(\lambda)$. By definition, we have 
$$e_c(\A_g,\VVl)=(-1)^{\frac{g(g+1)}{2}} \widehat{S}[n(\lambda)]+e_{g,{\rm endo}}(\lambda)+e_{g,\Eis}(\lambda).$$

We also define the {\sl extraneous} contribution $e_{g,\rm extr}(\lambda)$ through the following equation: 
$$e_c(\A_g,\VVl)=(-1)^{\frac{g(g+1)}{2}} S[n(\lambda)]+e_{g,{\rm extr}}(\lambda).$$

\subsection{Preview}
In Section~\ref{sec-g3} we will formulate a conjecture for the motivic Euler characteristic $ e_c(\A_3,\VVl)$ for any local system $\VVl$ on $\A_3$, in terms of $L$, $S[n_1]$, $S[n_1,n_2]$ and $S[n_1,n_2,n_3]$. This conjecture was found with the help of computer counts over finite fields: we have calculated
\begin{equation}\label{eq-trace} 
\Trace(F_q,e_c(\A_3 \otimes \overline{\QQ},{}_{\ell}\VVl)) 
\end{equation} 
for all prime powers $q \leq 17$ and all $\lambda$ with $|\lambda| \leq 60$.

In Section~\ref{sec-points}, we explain how we did these counts. In sections~\ref{sec-g3-charpol} and~\ref{sec-congruences}, we will discuss properties of Siegel modular cusp forms of degree three that we find under the assumption that our conjecture is true: we consider characteristic polynomials and the generalized Ramanujan conjecture as well as lifts from $\Gtwo$ and various congruences.

Before this, in Section~\ref{sec-g2}, we will review the situation in genus~$2$, with and without level $2$ structure, which we dealt with in the articles \cite{FvdG} and \cite{BFG}.

%%%%%%%%%%%%%%%%%%%
\section{Genus two} \label{sec-g2}
\subsection{The regular case} \label{sec-g2-r}
In the article \cite{FvdG}, the two latter authors formulated a conjectural analogue of Theorem~\ref{thm-g1} for genus~$2$ (as we will do here for genus~$3$ in Conjecture~\ref{conj-g3}). It was based on the integers
$$\Trace\bigl(F_q,e_c(\A_2 \otimes \overline{\QQ},{}_{\ell}\VVl)\bigr)$$
for all prime powers $q \leq 37$ and all $\lambda$ with $|\lambda| \leq 100$, which were found by counting points over finite fields, compare Section~\ref{sec-points}. Here is a reformulation of this conjecture, using the conjectural motive $\widehat{S}[n_1,n_2]$ described in Section~\ref{sec-motives}. 

\begin{conjecture} \label{conj-g2} The motivic Euler characteristic $e_c({\A}_2,\VVl)$ for regular $\lambda=(a,b)$ (with $a+b$ even) is given by
$$e_c({\A}_2,\VVl)=-\widehat{S}[a-b,b+3]+e_{2,\rm Eis}(a,b)+e_{2,\rm endo}(a,b),$$
where 
$$e_{2,{\rm Eis}}(a,b)=s_{a-b+2}-s_{a+b+4}L^{b+1}+\begin{cases} S[b+2] +1, 
& a \:\: \text{even}, \\ -S[a+3], & a \:\: \text{odd}, \end{cases}$$
and
$$e_{2,\rm endo}(a,b)= -s_{a+b+4} S[a-b+2] L^{b+1}.$$
Moreover, all regular $\lambda$ are normal. 
\end{conjecture}

The formula for $e_{2,\rm Eis}(a,b)$ has been proved by the third author in the categories $K_0(\MHSC)$ and $K_0(\GalC)$ using the BGG-complex, \cite[Corollary~9.2]{vdG2}. Moreover, Weissauer has proved the conjectured formula for the Euler characteristic of the inner cohomology in $K_0(\GalC)$, see \cite[Theorem 3]{Weissauer} (but note that in the formulation of the theorem, the factor $4$ should be removed). Consequently, Conjecture~\ref{conj-g2} is proved completely in $K_0(\GalC)$.

This gives us the possibility of computing traces of Hecke operators on spaces of Siegel modular forms of degree $2$ using counts of points over finite fields. 

\begin{example} For $\lambda=(11,7)$, the conjecture tells us that 
$$e_c({\A}_2,\VVl)=-S[4,10]-L^8.$$ 
We then know that for all primes $p \neq \ell$, 
$$\Trace(T(p),S_{4,10})=-\Trace(F_p,e_c\bigl(\A_2 \otimes \overline{\QQ},{}_{\ell}\VVl)\bigr)-p^8.$$
By Tsushima's dimension formula (see below), we have $s_{4,10}=1$. The counts of points over finite fields thus give the eigenvalues of any non-zero form $F$ in $S_{4,10}$ for all primes $p \leq 37$. We have for example: 
$$\lambda_2(F)= -1680, \qquad \lambda_3(F)= 55080, \qquad \lambda_{37}(F)=11555498201265580.$$
\end{example}

\subsection{The non-regular case} \label{sec-g2-nr} 
For the non-regular local systems, Conjecture~\ref{conj-g2} may fail. We refine the conjecture for a general local system in the following way. 

\begin{conjecture} \label{conj-g2-nr} The motivic Euler characteristic $e_c({\A}_2,\VVl)$ for any $\lambda=(a,b) \neq (0,0)$ (with $a+b$ even) is given by
$$e_c({\A}_2,\VVl)=-S[a-b,b+3]+e_{2,\rm extr}(a,b),$$
where 
\begin{multline*}
e_{2,\rm extr}(a,b) := -s_{a+b+4} S[a-b+2] L^{b+1}+\\
+s_{a-b+2}-s_{a+b+4}L^{b+1}+\begin{cases} S[b+2] +1, & a \:\: \text{even}, \\ -S[a+3], &  a \:\: \text{odd}. \end{cases}
\end{multline*}
\end{conjecture}

The difference between the general case and the regular case is that the Eisenstein and endoscopic contributions may behave irregularly. But we believe that their sum will not, except in the case when $\lambda$ is not normal, see Section~\ref{sec-g2-nonexh}.

There is a formula by Tsushima, see \cite{Tsushima}, for the dimension $s_{j,k}$ for all $j \geq 0$ and $k \geq 3$ which is proven for all $j \geq 1$ and $k \geq 5$ and for $j=0$ and $k \geq 4$ (note that the ring of scalar-valued modular forms and its ideal of cusp forms were determined by Igusa \cite{Igusa}). On the other hand, taking dimensions in Conjecture~\ref{conj-g2-nr}, we get that
\begin{multline} \label{eq-Euler-g2} 
-E_c({\A}_2,\VVl)-2 \, s_{a+b+4} \, s_{a-b+2}+{}\\{}+s_{a-b+2} -s_{a+b+4}+\begin{cases}2\, s_{b+2} +1, & a \:\: \text{even}, \\ -2\,s_{a+3}, & a \:\: \text{odd}, \end{cases}
\end{multline}
should equal $4 \, s_{a-b,b+3}$. By what was said in Section~\ref{sec-g2-r}, this is true for all regular $(a,b)$. 

We can decompose $\A_2 = \M_2 \cup \A_{1,1}$, where $\A_{1,1} := (\A_1\times \A_1)/\s_2$. There is a formula by Getzler for $E_c(\M_2,\VVl)$ for any $\lambda$, see \cite{Getz}. Together with the following formula, where $m:=(a-b)/2$ and $n:=(a+b)/2$,
\begin{multline*}
E_c(\A_{1,1},\VV_{a,b})=
\sum_{i=0}^b \, \sum_{j=0}^{m-1} E_c(\A_1,\VV_{a-i-j}) \, E_c(\A_1,\VV_{b-i+j}) +{}\\ +\sum_{k=m}^n \begin{cases} E_c(\A_1,\VV_{k}) \, \bigl(E_c(\A_1,\VV_{k})+1 \bigr)/2, & a \: \: \text{even}, \\ E_c(\A_1,\VV_{k}) \, \bigl(E_c(\A_1,\VV_{k})-1\bigr)/2, & a \: \: \text{odd}, \end{cases} 
\end{multline*}
we get a formula for $E_c(\A_2,\VVl)$. Grundh (see \cite{Grundh}) recently proved that for all $(a,b)$, Formula \eqref{eq-Euler-g2} is equivalent to Tsushima's dimension formula. From this, it follows that Tsushima's dimension formula also holds when $k=4$ and $j>0$.

\begin{example} For all $\lambda$ for which \eqref{eq-Euler-g2} tells us that $s_{a-b,b+3}=0$ (there are $85$ cases), we find as expected that 
$$\Trace\bigl(F_p,e_c({\A}_2,\VVl)\bigr)=\Trace\bigl(F_p,e_{2,\rm extr}(a,b)\bigr)$$
for all primes $p \leq 37$.
\end{example}

\begin{example} For $\lambda=(50,0)$, Conjecture~\ref{conj-g2-nr} tells us that 
$$e_c({\A}_2,\VVl)=-S[50, 3]-4 \, S[52] L+4-5 L.$$ 
For all primes $p \neq \ell$, it should then hold that 
$$\Trace(T(p),S_{50,3})=-\Trace(F_p,e_c\bigl(\A_2 \otimes \overline{\QQ},{}_{\ell}\VVl)\bigr)-4p\cdot \Trace(T(p),S_{52})+4-5p.$$ 
Formula~\eqref{eq-Euler-g2}, together with the fact that $E_c({\A}_2,\VV_{(50,0)})=-37$ (see \cite{Getz}), then tells us that (conjecturally) $s_{50,3}=1$. The counts of points over finite fields thus give the eigenvalues of any non-zero form $F$ in $S_{50,3}$ for all primes $p \leq 37$. We have for example (conjecturally): 
\begin{gather*} 
\lambda_2(F)= -37528320, \qquad \lambda_3(F)=-3184692509880, \\ \lambda_{37}(F)=86191557288628492956664825102803613165420.
\end{gather*}
\end{example}

\begin{example} For $\lambda=(32,32)$, Conjecture~\ref{conj-g2-nr} tells us that 
$$e_c({\A}_2,\VVl)=-S[0, 35]+S[34]+5 L^{34}$$
and the space $S_{0,35}$ is generated by one form, which Igusa denoted $\chi_{35}$. Our conjecture then tells us for instance that 
\begin{gather*}
\lambda_2(\chi_{35})= -25073418240, \qquad \lambda_3(\chi_{35})=-11824551571578840, \\ \lambda_{37}(\chi_{35})=-47788585641545948035267859493926208327050656971703460.
\end{gather*}
\end{example}

\subsection{Weight zero} \label{sec-wt0-g2}
As in the case of genus~$1$, we want Conjecture~\ref{conj-g2-nr} to hold also for $\QQl$-coefficients. Recall that $e_c({\A}_2,\VV_{0,0})=L^3+L^2$ and that these two classes belong to the Eisenstein cohomology. We thus extend the conjecture to the case $\lambda=(0,0)$ by defining $S[0,3]:=-L^3-L^2-L-1$ and $s_{0,3}=-1$ (remembering that $S[2]=-L-1$ and $s_2=-1$). Since the inner cohomology vanishes, we have $\widehat{S}[0,3]=0$. 

\subsection{Saito-Kurokawa lifts} \label{sec-g2-nonexh} 
For $\lambda$ regular, every Siegel cusp form of degree $2$ and weight $n(\lambda)$ that is a Hecke eigenform gives rise to a $4$-dimensional piece of the inner cohomology. We believe that this fails for the special kind of Siegel modular forms called Saito-Kurokawa lifts, which should contribute only $2$-dimensional pieces. In their presence, $\lambda$ will thus fail to be normal. 

For $a$ odd, the Saito-Kurokawa lift is a map $S_{2a+4} \to S_{0,a+3}$, see \cite{Z}. We can split $S_{0,a+3}$ as an orthogonal direct sum of the Maass Spezialschar $\Sigma^{SK}_{0,a+3}$ and its orthocomplement (w.r.t.~the Petersson inner product); we denote the latter space by $\Sigma_{0,a+3}^{\gen}$ and observe that it is stable under the Hecke algebra. The Saito-Kurokawa lift $F$ of a Hecke eigenform $f \in S_{2a+4}$ is a Hecke eigenform in $S_{0,a+3}$ with spinor $L$-function 
$$\zeta(s-a-1)\zeta(s-a-2)L(f,s).$$ 
This tells us that for all $p$, the trace of $T(p)$ on the Maass Spezialschar equals 
$$\Trace(T(p),S_{2a+4})+s_{2a+4}(p^{a+2}+p^{a+1}).$$

We conjecture that $\Sigma_{0,a+3}^{\gen}$ has the expected properties (see Section~\ref{sec-motives}), with $\Sigma^{\gen}[0,a+3]$ the corresponding piece of $\widehat{S}[0,a+3]$. The contribution corresponding to $\Sigma^{SK}_{0,a+3}$ will be denoted by $\Sigma^{SK}[0,a+3]$, and the first guess would be that $\Sigma^{SK}[0,a+3]$ equals $S[2a+4]+s_{2a+4}(L^{a+2}+L^{a+1})$. Considering the Hodge types, we find that $S[2a+4]$, $s_{2a+4} L^{a+2}$ and $s_{2a+4} L^{a+1}$ have to live inside $H^i_!(\A_2,\VV_{a,a})$ for $i=3,4,2$ respectively. This is not possible, because they should all have the same sign. Instead, we conjecture that $\Sigma^{SK}[0,a+3]=S[2a+4]$, compare Section~\ref{sec-g2-lev2} and \cite{Taylor-Siegel3}. 

\begin{conjecture} For $\lambda=(a,a)$, where $a$ is odd, we have 
$$\widehat{S}[0,a+3]= \Sigma^{SK}[0,a+3]+\Sigma^{\gen}[0,a+3] = S[2a+4]+\Sigma^{\gen}[0,a+3].$$
In all other cases, $\lambda$ is normal. 
\end{conjecture}

It follows from this conjecture that for $a$ odd,
$$\Trace(F_p,\widehat{S}[0,a+3])=\Trace(T(p),S_{0,a+3})-s_{2a+4}(p^{a+2}+p^{a+1}).$$

\begin{example} For $\lambda=(11,11)$, Conjecture~\ref{conj-g2-nr} tells us that 
$$e_c({\A}_2,\VVl)=-S[0,14]-1+L^{13}.$$ 
Formula~\eqref{eq-Euler-g2}, together with the fact that $E_c({\A}_2,\VV_{(11,11)})=-4$ (see \cite{Getz}), then tells us that (conjecturally) $s_{0,14}=1$. Since $s_{26}=1$, the modular form in $S_{0,14}$ is a Saito-Kurokawa lift. We should then have that 
$$\Trace(F_p,e_c({\A}_2,\VVl))=-\bigl(\Trace(T(p),S[26])+p^{12}+p^{13}\bigr)-1+p^{13}.$$ 
This indeed holds for all primes $p \leq 37$. 

For $a\leq 15$ odd, the space $S_{0,a+3}$ is generated by Saito-Kurokawa lifts. In all these cases $\Trace(F_p,e_c({\A}_2,\VV_{a,a}))$ equals 
$$-\bigl(\Trace(F_p,S[2a+4])+s_{2a+4}(p^{a+1}+p^{a+2})\bigr)+\Trace\bigl(F_p,e_{2,\rm extr}(a,a)\bigr)$$ 
for all primes $p \leq 37$. 
\end{example}

\subsection{Characteristic polynomials} \label{sec-g2-charpol} 
We would like to state some observations on the factorization of the local spinor $L$-factor of Siegel modular forms of degree $2$. 

Any Siegel modular form of degree $2$ that is not a Saito-Kurokawa lift gives rise to a four-dimensional piece of the corresponding third inner cohomology group (see \cite{Weissauer} and \cite{Taylor-Siegel3}). We will disregard the Saito-Kurokawa lifts (since their spinor $L$-functions will come from elliptic modular forms) and consider the spinor $L$-functions of elements in $\Sigma^{\gen}_{a-b,b+3}$. For an eigenform $f \in \Sigma^{\gen}_{a-b,b+3}$, the local spinor $L$-factor $L_p(s,f)$ equals $Q_p(p^{-s},f)^{-1}$, where $Q_p(X,f)$ is the following polynomial: $$1-\lambda_p(f)\,X + \bigl(p\lambda_{1,p^2}(f)+(p^3+p) \lambda_{2,p^2}(f)\bigr)X^2 - \lambda_p(f)\,p^{a+b+3} \, X^3 + p^{2(a+b+3)}X^4.$$ 
All forms in $\Sigma^{\gen}_{a-b,b+3}$ will fulfil the Ramanujan conjecture, that is, for every prime $p \neq \ell$, the roots of the characteristic polynomial just stated will have absolute value $p^{-(a+b+3)/2}$, see \cite[Th.~3.3.3]{Weissauer-book}.

The polynomial $Q_p(X,f)$ for $f \in \Sigma^{\gen}_{a-b,b+3}$ is equal to the characteristic polynomial of Frobenius acting on the corresponding $\ell$-adic representations, which are found inside $\Sigma^{\gen}[a-b,b+3]$. To compute the characteristic polynomial of Frobenius $F_p$ acting on $\Sigma^{\gen}[a-b,b+3]$, we therefore need to compute $\Trace(F^i_p,\Sigma^{\gen}[a-b,b+3])$ for $1 \leq i \leq 2\, \dim \Sigma^{\gen}_{a-b,b+3}$. Using Conjecture~\ref{conj-g2-nr} and the fact that $F^i_p=F_{p^i}$, we can reformulate this to computing $\Trace(F_{p^i},e_c({\A}_2,\VV_{a,b}))$ for $i$ from $1$ to $2\, \dim \Sigma^{\gen}_{a-b,b+3}$. 

Since we have computed $\Trace(F_{2^i},e_c({\A}_2,\VV_{a,b}))$ for $1 \leq i \leq 4$, we can determine the characteristic polynomial at $p=2$ for any $\lambda=(a,b)$ such that $\dim \Sigma^{\gen}_{a-b,b+3} \leq 2$. For $|\lambda| \leq 100$, there are $40$ choices of $\lambda=(a,b)$ for which $s^{\gen}_{a-b,b+3}:=\dim \Sigma^{\gen}_{a-b,b+3}=1$ and $27$ choices for which $s^{\gen}_{a-b,b+3}= 2$ (this follows from the formula in Section~\ref{sec-g2-nr}). These characteristic polynomials are irreducible over $\QQ$, except for the following local systems: $\lambda=(22,4),(20,10),(21,21),(23,23)$. The factorization into two polynomials of degree $4$ for the two latter local systems was found by Skoruppa \cite{Skoruppa}. It is shown in \cite{Ram-Sha} that there is a Siegel modular form $f$ of degree $2$ and weight $(10,13)$ such that $\lambda_p(f)=\Trace\bigl(F_p,\Sym^3(S[12])\bigr)$ for all primes $p$. This accounts for the splitting of the characteristic polynomial in the case $\lambda=(20,10)$. In the case $\lambda=(22,4)$, the characteristic polynomial at $p=2$ splits in the following way:
\begin{multline*} 
(1+ 32736 \, X + 857571328 \, X^2+ 32736 \cdot 2^{29} \, X^3+ 2^{58} \, X^4) \cdot \\ (1- 7920 \, X+45752320  \, X^2-7920  \cdot 2^{29} \, X^3+ 2^{58} \, X^4).
\end{multline*}

\subsection{The case of level $2$} \label{sec-g2-lev2} 
In \cite{BFG}, a similar kind of analysis was made for local systems on the moduli space ${\mathcal A}_2[2]$ of principally polarized abelian surfaces together with a full level two-structure. Here one has the additional structure of the action of $\GSp(4,\ZZ/2) \cong \s_6$ on the cohomology groups. We will give some additional comments on the (conjectural) picture we have in this case. Our understanding benefited very much from two letters sent to us by Deligne (\cite{Deligne-letters}). 

What was called middle endoscopy in~\cite{BFG} is called endoscopy in this article. For any two Hecke eigenforms of level $2$, $f \in S_{a+b+4}(\Gamma(2))$ and $g \in S_{a-b+2}(\Gamma(2))$, there is a contribution to $e_c({\A}_2[2],\VV_{a,b})$ of the form  
$$X \otimes M_f + Y \otimes L^{b+1} \, M_g,$$
where $X$ and $Y$ are (either zero or) different representations of $\s_6$. The contributions of the form $Y \otimes L^{b+1} \, M_g$ make up the endoscopic cohomology. The contributions of the form $X \otimes M_f$ correspond to the cases where there is a lift of Yoshida type (see \cite[Conj.~6.1, 6.4]{BFG}), taking the two forms $f$ and $g$ to an eigenform $F \in S_{a-b,b+3}(\Gamma_2)$ with spinor $L$-function 
$$L(F,s)=L(f,s) \, L(g,s-b-1).$$
In the cohomology, we thus only see the two-dimensional piece $M_f$, corresponding to the factor $L(f,s)$ of the $L$-function, instead of the expected four-dimensional piece (compare Section~\ref{sec-g2-nonexh}). Note that these non-exhaustive lifts in level $2$ occur for regular local systems. 

If $a=b$, we can take $M_g$ equal to $L+1$. The corresponding contributions are then of Saito-Kurokawa type (cf.~\cite[Conj.~6.1, 6.6]{BFG}).

\begin{example} In the notation of \cite{BFG}, $e_c({\A}_2[2],\VVl)$ for $\lambda=(5,1)$ equals (conjecturally) the sum of the following contributions. Here $s[\mu]$ denotes the irreducible representation of $\s_6$ corresponding to the partition $\mu$. First there is the Eisenstein cohomology, 
\begin{multline*}
-S[\Gamma_0(2),8]^{\mathrm{new}}(s[2^3]+s[3,2,1]+s[4,2])\\-L^2 (s[3,2,1]+s[3^2]+s[4,1^2]+s[4,2]+s[5,1])+(s[3^2]+s[4,1^2]),
\end{multline*}
then the endoscopic cohomology, 
$$-L^2 \, S[\Gamma_0(4),6]^{\mathrm{new}}(s[3,1^3]+s[4,1^2]),$$
and finally there is a lift of Yoshida type, contributing 
$$-S[\Gamma_0(4),10]^{\mathrm{new}} \, s[2,1^4].$$
\end{example}

%%%%%%%%%%%%%%%%%%%%%
\section{Genus three} \label{sec-g3} 
We now formulate our main conjecture. 

\begin{conjecture} \label{conj-g3} The motivic Euler characteristic $e_c({\A}_3,\VVl)$ for any $\lambda =(a,b,c) \neq (0,0,0)$ is given by
$$e_c({\A}_3,\VV_{a,b,c})=S[a-b,b-c,c+4]+e_{3,\rm extr}(a,b,c),$$
where 
\begin{multline*} 
 e_{3,\rm extr}(a,b,c):= -e_c({\A}_2,\VV_{a+1,b+1})+e_c({\A}_2,\VV_{a+1,c})-e_c({\A}_2,\VV_{b,c}) \\
-e_{2,\rm extr}(a+1,b+1) \otimes  S[c+2] + e_{2,\rm extr}(a+1,c) \otimes S[b+3] \\-e_{2,\rm extr}(b,c)\otimes S[a+4].
\end{multline*}
\end{conjecture}

\begin{remark} The term $e_{3,\rm extr}(a,b,c)$ is thus formulated in terms of genus~$2$ contributions, which can be computed using the results of Section~\ref{sec-g2}. The Euler characteristics for $\lambda$ with $|\lambda|\leq18$ are given at the end of the paper.
\end{remark}

\subsection{Integer Euler characteristics}
The numerical Euler characteristic $E_c({\A}_3,\VVl)$, which was defined in Section~\ref{sec-localsystems}, can be computed for any $\lambda$ using the results of \cite{BvdG,JBvdG}. Taking dimensions in the formula in Conjecture~\ref{conj-g3} we end up with the following definition:
\begin{multline*}
 E_{3,\rm extr}(a,b,c):=-E_c({\A}_2,\VV_{a+1,b+1})+E_c({\A}_2,\VV_{a+1,c})-E_c({\A}_2,\VV_{b,c}) \\
- 2 \, E_{2,\rm extr}(a+1,b+1) \, s_{c+2}  +  2\,  E_{2,\rm extr}(a+1,c)\, s_{b+3}  -2 \, E_{2,\rm extr}(b,c)\, s_{a+4}\,, 
\end{multline*}
and the following conjecture.
\begin{conjecture} \label{conj-g3-int} For any local system $\lambda=(a,b,c)\neq (0,0,0)$, 
$$s_{a-b,b-c,c+4}=\frac{1}{8} \Bigl(E_c({\A}_3,\VV_{a,b,c})-E_{3,\rm extr}(a,b,c)\Bigr).$$ 
\end{conjecture}

The first check of this conjecture is that $E_c({\A}_3,\VVl)-E_{3,\rm extr}(\lambda)$ is a nonnegative integer divisible by $8$ for each $\lambda$ with $|\lambda| \leq 60$. The dimension of the space $S_{0,0,k}$ of classical (i.e., scalar-valued) Siegel modular cusp forms of weight $k$ is known by Tsuyumine (for $k\geq 4$), see \cite{Tsuyumine}. In the case $\lambda=(a,a,a)$, we have then checked that Conjecture~\ref{conj-g3-int} is true for all $0 \leq a \leq 20$. 

\begin{table}[ht]\caption{Conjectural dimensions of spaces of degree $3$ cusp forms.} \label{tab-eulg3} 
\vbox{
\centerline{\def\quad{\hskip 0.3em\relax}
\vbox{\offinterlineskip
\hrule
\halign{&\vrule#& \quad \hfil#\hfil \strut \quad \cr
height2pt&\omit&&\omit&&\omit&&\omit&&\omit&&\omit&&\omit&&\omit&\cr
& $\lambda$ && $s_{n(\lambda)}$ && $\lambda$ && $s_{n(\lambda)}$ && $\lambda$ && $s_{n(\lambda)}$ && $\lambda$ && $s_{n(\lambda)}$ &\cr
height2pt&\omit&&\omit&&\omit&&\omit&&\omit&&\omit&&\omit&&\omit&\cr
\noalign{\hrule}
height2pt&\omit&&\omit&&\omit&&\omit&&\omit&&\omit&&\omit&&\omit&\cr
& $(20, 0, 0)$ && $0$ && $(19, 1, 0)$ && $0$ && $(18, 2, 0)$ && $0$ && $(18, 1, 1)$ && $0$ &\cr
& $(17, 3, 0)$ && $0$ && $(17, 2, 1)$ && $0$ && $(16, 4, 0)$ && $0$ && $(16, 3, 1)$ && $0$ &\cr
& $(16, 2, 2)$ && $0$ && $(15, 5, 0)$ && $0$ && $(15, 4, 1)$ && $0$ && $(15, 3, 2)$ && $0$ &\cr
& $(14, 6, 0)$ && $0$ && $(14, 5, 1)$ && $0$ && $(14, 4, 2)$ && $1$ && $(14, 3, 3)$ && $0$ &\cr
& $(13, 7, 0)$ && $0$ && $(13, 6, 1)$ && $1$ && $(13, 5, 2)$ && $0$ && $(13, 4, 3)$ && $1$ &\cr
& $(12, 8, 0)$ && $0$ && $(12, 7, 1)$ && $0$ && $(12, 6, 2)$ && $1$ && $(12, 5, 3)$ && $0$ &\cr
& $(12, 4, 4)$ && $1$ && $(11, 9, 0)$ && $0$ && $(11, 8, 1)$ && $0$ && $(11, 7, 2)$ && $0$ &\cr
& $(11, 6, 3)$ && $1$ && $(11, 5, 4)$ && $0$ && $(10, 10, 0)$&& $0$ && $(10, 9, 1)$ && $1$ &\cr
& $(10, 8, 2)$ && $1$ && $(10, 7, 3)$ && $0$ && $(10, 6, 4)$ && $1$ && $(10, 5, 5)$ && $0$ &\cr
& $(9, 9, 2)$  && $0$ && $(9, 8, 3)$  && $0$ && $(9, 7, 4)$  && $0$ && $(9, 6, 5)$  && $0$ &\cr
& $(8, 8, 4)$  && $0$ && $(8, 7, 5)$  && $0$ && $(8, 6, 6)$  && $1$ && $(7, 7, 6)$  && $0$ &\cr
& $(22, 0, 0)$ && $0$ && $(21, 1, 0)$ && $0$ && $(20, 2, 0)$ && $0$ && $(20, 1, 1)$ && $0$ &\cr
& $(19, 3, 0)$ && $0$ && $(19, 2, 1)$ && $0$ && $(18, 4, 0)$ && $0$ && $(18, 3, 1)$ && $0$ &\cr
& $(18, 2, 2)$ && $0$ && $(17, 5, 0)$ && $0$ && $(17, 4, 1)$ && $1$ && $(17, 3, 2)$ && $0$ &\cr
& $(16, 6, 0)$ && $0$ && $(16, 5, 1)$ && $1$ && $(16, 4, 2)$ && $1$ && $(16, 3, 3)$ && $0$ &\cr
& $(15, 7, 0)$ && $0$ && $(15, 6, 1)$ && $1$ && $(15, 5, 2)$ && $1$ && $(15, 4, 3)$ && $1$ &\cr
& $(14, 8, 0)$ && $0$ && $(14, 7, 1)$ && $1$ && $(14, 6, 2)$ && $1$ && $(14, 5, 3)$ && $1$ &\cr
& $(14, 4, 4)$ && $1$ && $(13, 9, 0)$ && $1$ && $(13, 8, 1)$ && $1$ && $(13, 7, 2)$ && $1$ &\cr
& $(13, 6, 3)$ && $1$ && $(13, 5, 4)$ && $1$ && $(12, 10, 0)$&& $0$ && $(12, 9, 1)$ && $0$ &\cr
& $(12, 8, 2)$ && $1$ && $(12, 7, 3)$ && $1$ && $(12, 6, 4)$ && $1$ && $(12, 5, 5)$ && $0$ &\cr
& $(11, 11, 0)$&& $0$ && $(11, 10, 1)$&& $1$ && $(11, 9, 2)$ && $1$ && $(11, 8, 3)$ && $1$ &\cr
& $(11, 7, 4)$ && $1$ && $(11, 6,5)$  && $0$ && $(10, 10, 2)$&& $0$ && $(10, 9, 3)$ && $0$ &\cr
& $(10, 8, 4)$ && $1$ && $(10, 7, 5)$ && $0$ && $(10, 6, 6)$ && $0$ && $(9, 9,4)$   && $0$ &\cr
& $(9, 8, 5)$  && $0$ && $(9, 7, 6)$  && $0$ && $(8, 8, 6)$  && $0$ && $(8, 7, 7)$  && $0$ &\cr
} \hrule}
}}
\end{table}

\begin{example} There are $317$ choices of $\lambda$ for which $|\lambda| \leq 60$ and Conjecture~\ref{conj-g3-int} tells us that $s_{n(\lambda)}=0$. For all these choices, we find that $\Trace \bigl(F_q,e_c({\A}_3,\VVl)\bigr)=\Trace \bigl(F_q,e_{3,\rm extr}(\lambda)\bigr)$ for all $q \leq 17$.

One such instance is $\lambda=(15,3,0)$, where Conjecture~\ref{conj-g3} tells us that 
$$e_{3,\rm extr}(\lambda)=S[12,7]-S[18]L-2L^6-L+1.$$
\end{example}

\begin{example} For $|\lambda| \leq 18$, the only cases for which Conjecture~\ref{conj-g3-int} tells us that $s_{n(\lambda)}$ is non-zero are $\lambda=(8,4,4)$, $(11,5,2)$ and $(9,6,3)$, and then $s_{n(\lambda)}=1$. In Table~\ref{tab-eulg3} we list $s_{n(\lambda)}$ for all $\lambda$ of weight $20$ and $22$.

In the case $\lambda=(11,5,2)$, Conjecture~\ref{conj-g3-int} tells us that 
$$e_c({\A}_3,\VVl)=S[6,3,6]-S[12] \, L^3+L^7-L^3+1.$$
For all primes $p \neq \ell$ we should then have that 
$$\Trace(T(p),S_{6,3,6})=\Trace\bigl(F_p,e_c({\A}_3,\VVl)\bigr)-\bigl(-\Trace(T(p),S_{12}) \, p^3+p^7-p^3+1\bigr).$$
The space $S_{6,3,6}$ should be spanned by one form, say $F_{6,3,6}$, and from our computations we get that (conjecturally)
$$\lambda_2(F_{6,3,6})=0, \quad \lambda_3(F_{6,3,6})=-453600, \quad  \lambda_{17}(F_{6,3,6})=-107529004510200.$$
\end{example}

\begin{example} For $\lambda=(15,13,12)$, Conjecture~\ref{conj-g3-int} tells us that $s_{n(\lambda)}=1$ and Conjecture~\ref{conj-g3-int} that 
$$e_c({\A}_3,\VVl)=S[2,1,16]-S[4,15]-2S[16]\,L^{13}+2L^{15}-2L^{13}.$$
Conjecturally, for $F_{2,1,16}$ a generator of $S_{2,1,16}\,$, we find that 
\begin{gather*}
  \lambda_2(F_{2,1,16})=6994944, \qquad \lambda_3(F_{2,1,16})=134431309152, \\  \lambda_{17}(F_{2,1,16})=14399876302755866405698174344.
\end{gather*}
\end{example}

\subsection{Weight zero} \label{sec-wt0-g3}
In order to make Conjecture~\ref{conj-g3} true for $\lambda=(0,0,0)$ we put 
$$S[0,0,4]:=L^6+L^5+L^4+2\, L^3+L^2+L+1, $$
because we know that $e_c({\A}_3,\VV_{0,0,0})=L^6+L^5+L^4+L^3+1$, see \cite{Hain}. All of this cohomology is Eisenstein and thus $\widehat{S}[0,0,4]=0$. 

Note that by our definitions, the rank of $S[2]$ equals $-2$, the rank of $S[0,3]$ equals $-4$, and the rank of $S[0,0,4]$ equals $8$. The definitions so far fit into the pattern
$$S[0,0,\ldots,0,g+1]=(-1)^g (1+L^g) \, S[0,0,\ldots,0,g],$$
but this will not be true in general, because at some point classical cusp forms of weight $g+1$ will appear. In fact, it follows from work of Ikeda~\cite{Ikeda} that there exists a classical cusp form of weight $12$ in genus~$11$ (e.g., the form $F_{23}$ in Table 3, p.~494).

\subsection{Non-exhaustive lifts} \label{sec-g3-nonexh}
We conjecture the existence of the following three types of lifts of Siegel modular forms. 

\begin{conjecture}\label{conj-g3-lifts} The Hecke module $S_{a-b,b-c,c+4}$ splits into a direct sum of the two submodules $\Sigma_{a-b,b-c,c+4}^{\gen}$ and $\Sigma_{a-b,b-c,c+4}^{\rm ne}\,$, where the latter is generated by the following forms. 
\item{i)} For any choice of eigenforms $f \in S_{b+3}$, $g \in S_{a+c+5}$ and $h\in S_{a-c+3}$ there exists an eigenform $F \in S_{a-b,b-c,c+4}$ with spinor $L$-function 
$$L(F,s)=L(f\otimes g, s)L(f\otimes h, s-c-1).$$
\item{ii)} For any choice of eigenforms $f \in S_{a+4}$ and $g \in S_{2b+4}$ there exists an eigenform $F \in S_{a-b,0,b+4}$ with spinor $L$-function 
$$L(F,s)=L(f,s-b-1)L(f,s-b-2)L(f\otimes g,s).$$
\item{iii)} For any choice of eigenforms $f \in S_{c+2}$ and $g \in S_{2a+6}$ there exists an eigenform $F \in S_{0,a-c,c+4}$ with spinor $L$-function 
$$L(F,s)=L(f,s-a-2)L(f,s-a-3)L(f\otimes g,s).$$
\end{conjecture}

A precursor of this conjecture for the case of classical Siegel modular forms can be found in the work of Miyawaki. He observed in \cite{Miyawaki} that for the non-zero cusp form $F_{12}$ of weight $(0,0,12)$, the Euler factor at $p=2$ of the spinor $L$-function $L(F_{12},s)$ was equal to $L_2(\Delta,s-9) L_2(\Delta,s-10) L_2(\Delta\otimes g,s)$, with eigenforms $\Delta \in S_{12}$ and $g \in S_{20}$, see Example~\ref{exa-miyawaki}. He then conjectured the equality of spinor $L$-functions 
\begin{equation} \label{eq-miyawaki}
L(F_{12},s)=L(\Delta,s-9) L(\Delta,s-10)L(\Delta\otimes g,s).
\end{equation}
This conjecture was proved by Ikeda \cite{Ikeda}, p.~474. 

Miyawaki made a similar observation for the cusp form $F_{14}$ of weight $(0,0,14)$, where now the role of $g$ is taken by an eigenform $h \in S_{26}$. He also observed that such a lift of a pair $f_1,f_2$ of cusp forms of weight $k_1$ and $k_2$ to scalar-valued Siegel modular forms of weight $(0,0,k)$ could only occur if $(k_1,k_2)$ was equal to $(k,2k-4)$ or $(k-2,2k-2)$. He then conjectured part (ii) and (iii) of Conjecture~\ref{conj-g3-lifts} for all $k-4=a=b=c$. 

\begin{example} Let us check the two examples found by Miyawaki. In the case $\lambda=(8,8,8)$, Conjecture~\ref{conj-g3} tells us that
$$e_{3,\rm extr}(\lambda)=-L^{10}\, S[12]+S[0,12]-L^{11}-L^{10}+S[12]+1$$
and we can check that
$$\Trace(F_p,e_c({\A}_3,\VVl))-\Trace(F_p,e_{3,\rm extr}(\lambda))=\Trace\bigl(F_p,S[12] \otimes (S[20]+L^{10}+L^9)\bigr)$$
holds for all $p \leq 17$. Similarly, for $\lambda=(10,10,10)$ we have
$$e_{3,\rm extr}(\lambda)=-S[12] \, L^{13}+S[0,14]-L^{13}-2 \, L^{12}+1$$
and
$$\Trace\bigl(F_p,e_c({\A}_3,\VVl)\bigr)-\Trace\bigl(F_p,e_{3,\rm extr}(\lambda)\bigr)=\Trace\bigl(F_p,S[12] \otimes (S[26]+L^{13}+L^{12})\bigr)$$
holds for all $p \leq 17$.

Using Conjecture~\ref{conj-g3-int} and Conjecture~\ref{conj-g3-lifts}, we find $19$ choices, presented in Table~\ref{tab-lifts}, of $\lambda$ for which $|\lambda| \leq 60$ and $S_{n(\lambda)}=\Sigma_{n(\lambda)}^{\rm ne}$. In all these cases the expected equalities hold, just as for $\lambda=(8,8,8)$ and $(10,10,10)$. 
\end{example}

\begin{table}[ht] \caption{Local systems for which $S_{n(\lambda)}=\Sigma_{n(\lambda)}^{\rm ne}$.}\label{tab-lifts} 
\vbox{
\centerline{\def\quad{\hskip 0.3em\relax}
\vbox{\offinterlineskip
\hrule
\halign{&\vrule#& \quad \hfil#\hfil \strut \quad \cr
height2pt&\omit&&\omit&\cr
& $(a,b,c)$ && $S[a-b,b-c,c+4]$ &\cr 
\noalign{\hrule}
height2pt&\omit&&\omit&\cr
& $(8,4,4)$ && $S[12](S[12]+L^6+L^5)$ &\cr
& $(12,4,4)$ && $S[16](S[12]+L^6+L^5)$ & \cr
& $(10,9,1)$ && $S[12](S[16]+L^2S[12])$ &\cr
& $(8,6,6)$ && $S[12](S[16]+L^8+L^7)$ & \cr
& $(14,4,4)$ && $S[18](S[12]+L^6+L^5 $ & \cr
& $(13,9,0)$ && $S[12](S[18]+LS[16])$ & \cr
& $(11,9,2)$ && $S[12](S[18]+L^3S[12])$ & \cr
& $(15,9,0)$ && $S[12](S[20]+LS[18])$ & \cr
& $(8,8,8)$ && $S[12](S[20]+L^{10}+L^9)$ & \cr
& $(13,13,0)$ && $S[16](S[18]+LS[16])$ &\cr 
& $(15,13,0)$ && $S[16](S[20]+LS[18])$ & \cr
& $(13,13,4)$ && $S[16](S[22]+L^5S[12])$ &\cr 
& $(10,10,10)$ && $S[12](S[26]+L^{13}+L^{12})$ & \cr
& $(15,15,2)$ && $S[18](S[22]+L^3S[16])$ &\cr 
& $(12,10,10)$ && $S[16](S[24]+2(L^{12}+L^{11}))$ & \cr
& $(12,12,10)$ && $S[12](S[30]+2(L^{15}+L^{14}) $ &\cr 
& $(18,9,9)$ && $S[12](S[32]+2L^{10}S[12])$ & \cr
& $(14,14,14)$ && $S[16](S[34]+2(L^{17}+L^{16}))$ &\cr
& && $+S[18](S[32]+2(L^{16}+L^{15}))$ & \cr
 } \hrule}
}}
\end{table}

\begin{example} In the case $\lambda=(12,12,12)$, Conjecture~\ref{conj-g3-int} tells us that $s_{0,0,16}=3$ and the space generated by the lifts given in Conjecture~\ref{conj-g3-lifts} is two-dimensional. There is therefore a Hecke eigenform $F_{0,0,16}$ in $\Sigma^{\gen}_{0,0,16}$ whose eigenvalue for $T(p)$ is found (conjecturally) by the formula
$$\Trace\bigl(F_p,e_c({\A}_3,\VVl)\bigr)-\Trace\bigl(F_p,e_{3,\rm extr}(\lambda)\bigr)-\Trace\bigl(F_p,S[16] \otimes (S[28]+2\, L^{14}+2 \, L^{13})\bigr),$$
and we find for example:
\begin{gather*}
\lambda_2(F_{0,0,16})=-115200, \qquad \lambda_3(F_{0,0,16})=14457333600, \\ \lambda_{17}(F_{0,0,16})=-84643992509680105660020600.
\end{gather*}
\end{example}

\begin{remark} There are $123$ choices of $\lambda$ such that $|\lambda| \leq 60$ and for which the conjectures tell us that $s^{\gen}_{\lambda}:=\dim \Sigma_{n(\lambda)}^{\gen}=1$. 
\end{remark}

Consider an element of the Grothendieck group of motives of the form
$$M_f \otimes M_g + L^r \, M_f \otimes M_h,$$
where $f \in S_k$, $g \in S_l$ and $h \in S_m$ are eigenforms, or $m=2$ and $M_h=L+1$. Fix $\lambda=(a,b,c)$ and assume that $(k,l,m,r)$ is equal to one of the following three possibilities:
\begin{itemize}
\item[(i)] $(b+3,a+c+5,a-c+3,c+1)$; 
\item[(ii)] $(a+4,b+c+4,b-c+2,c+1)$;
\item[(iii)] $(c+2,a+b+6,a-b+2,b+2)$.
\end{itemize}
Assume furthermore that $(k,l,m,r)$ is such that there is a lift as in Conjecture~\ref{conj-g3-lifts}. One could then expect to find the whole $8$-dimensional motive in the cohomology, but, similarly to the genus~$2$ case, we conjecture that only the part $M_f \otimes M_g$ contributes. If $(k,l,m,r)$ is such that there is no lift as in Conjecture~\ref{conj-g3-lifts}, then we conjecture (in the regular case) that the terms $L^r \, M_f \otimes M_h$ will contribute to (and also generate) the endoscopic cohomology, see Section~\ref{sec-g3-regular}.

\begin{conjecture} \label{conj-g3-nonexh} For any local system $\lambda=(a,b,c)$, we have 
$$\widehat{S}[a-b,b-c,c+3]=\Sigma^{\gen}[a-b,b-c,c+4]+\Sigma^{\rm ne}[a-b,b-c,c+4]$$
where $\Sigma^{\gen}[a-b,b-c,c+4]$ has the expected properties (see Section~\ref{sec-motives}) and where $\Sigma^{\rm ne}[a-b,b-c,c+4]$ is the sum of the following three contributions:
\begin{itemize}
\item[(i)] $s_{a-c+3} \, S[b+3] \otimes S[a+c+5];$
\item[(ii)] $S[a+4] \otimes S[2b+4]$,  if $b=c\,;$
\item[(iii)] $S[c+2] \otimes S[2a+6]$,  if $a=b$ and $c>0$. 
\end{itemize}
\end{conjecture}

\subsection{The regular case} \label{sec-g3-regular}
In \cite{vdG2} it is proved that if $\lambda$ is regular, then the rank~$1$ part of the Eisenstein cohomology $e_{3,\Eis}(a,b,c)$ equals 
$$-e_c({\A}_2,\VV_{a+1,b+1})+e_c({\A}_2,\VV_{a+1,c})-e_c({\A}_2,\VV_{b,c}).$$ This is the first piece of the formula in Conjecture~\ref{conj-g3} for $e_{3,\rm extr}(a,b,c)$.

Let us make a refinement of Conjecture~\ref{conj-g3} in the regular case. Note that contributions to the endoscopic cohomology should have Deligne weight $a+b+c+6$ and should appear with positive sign. 

\begin{conjecture} \label{conj-g3-reg} If $\lambda=(a,b,c)$ is regular, then $e_{3,\rm endo}(a,b,c)$ is given by  
$$s_{b+c+4} \, S[a+4] S[b-c+2]L^{c+1}+s_{a+b+6} \, S[c+2] S[a-b+2] L^{b+2},$$
and $e_{3,\rm Eis}(a,b,c)$ by 
$$e_{3,\rm extr}(a,b,c)-e_{3,\rm endo}(a,b,c)+s_{a+c+5} \, S[b+3] S[a-c+3] L^{c+1}.$$
\end{conjecture}

\begin{example} For $\lambda=(16,13,3)$, Conjectures \ref{conj-g3-nonexh} and \ref{conj-g3-reg} read:
\begin{gather*}
\widehat{S}[3,10,7]=\Sigma^{\gen}[3,10,7]+S[16] \otimes S[24],\\
e_{3,\rm endo}(a,b,c)=S[20]\otimes S[12] \, L^4,\\
e_{3,\rm Eis}(a,b,c)=S[20] L^4-4S[16]L^4-2S[20]+S[12]L^4+2 \, S[16]-L^4.   
\end{gather*}
\end{example}

%%%%%%%%%%%%%%%%%%%%%%%%%%%%%%%%%%%%%%%%%%%%
\section{Counting points over finite fields} \label{sec-points}
In this section, we will indicate how we found the information necessary to compute the expression \eqref{eq-trace} for all $q \leq 17$ and $\lambda$ with weight  $|\lambda| \leq 60$. 

Just as in Section~\ref{sec-g1},  the trace \eqref{eq-trace} can be computed in terms of finite fields, that is, 
$$\Trace \bigl(F_q,e_c(\A_g \otimes \overline{\QQ},\VVl)\bigr) = \Trace \bigl(F_q,e_c(\A_g \otimes \overline{\QQ}_p,\VVl)\bigr) = \Trace\bigl(F_q,e_c(\A_g \otimes \FFb_q,\VVl)\bigr).$$

Let $[\A_g(\FF_q)]$ denote the set of $\FF_q$-isomorphism classes of principally polarized abelian varieties of dimension $g$ defined over $\FF_q$. For an $[A] \in [\A_g(\FF_q)]$, we denote by $\alpha_1(A),\ldots, \alpha_{2g}(A)$ the eigenvalues of the Frobenius map acting on $H^1(A\otimes \FFb_q,\QQl)$, ordered in such a way that $\alpha_i(A)\alpha_{g+i}(A)=q$. Let $s_{<\lambda>}(x_1,\ldots,x_g;t) \in \ZZ[t]$ be the Schur polynomial for $\Sp(2g,\QQ)$ associated to $\lambda$ and homogenized using $t$, where $x_1,\ldots,x_g$ have weight $1$ and $t$ weight~$2$, see \cite[p.~466]{FH}. We then find that 
$$\Trace(F_q,(\VVl)_{A\otimes \FFb_q})=s_{<\lambda>}\bigl(\alpha_1(A),\ldots,\alpha_{g}(A);q\bigr).$$
From the Lefschetz trace formula, it follows that (see \cite[Th.~3.2]{SGA})
\begin{equation} \label{eq-traceFq}
\Trace\bigl(F_q,e_c(\A_g \otimes \FFb_q,\VVl)\bigr) = \sum_{[A] \in [\A_g(\FF_q)]} \frac{s_{<\lambda>}\bigl(\alpha_1(A),\ldots,\alpha_{g}(A);q\bigr)}{|\Aut_{\FF_q}(A)|}.
\end{equation}

The monomials in the power sums, $p_k(x_1,\ldots,x_g):=\sum_i x_i^k$ for $1 \leq k \leq g$, form a rational basis for the symmetric polynomials in $x_1,\ldots,x_g$. For a partition $\mu$, let $\hat \mu$ denote the dual partition. We can then write 
$$s_{<\lambda>}=\sum_{|\mu| \leq |\lambda|} r^{\lambda}_{\mu} \,\, t^{\frac{|\lambda|-|\mu|}{2}} \,\,p_{\mu},
$$
for some $r^{\lambda}_{\mu} \in \QQ$ and where $p_{\mu}:=\prod_{i=1}^{\mu_1} p_{\hat \mu_i}$. By the Lefschetz trace formula, 
$$a_i(A):=\Trace\bigl(F_{q^i},H^1(A \otimes \FFb_q,\QQl)\bigr)=p_i \bigl(\alpha_1(A),\ldots,\alpha_{g}(A)\bigr)$$
for any $[A] \in [\A_g(\FF_q)]$, and so 
$$a_{\mu}(A):=\prod_{i=1}^{\mu_1} \Trace\bigl(F_{q^{\hat \mu_i}},H^1(A \otimes \FFb_q,\QQl)\bigr)=p_{\mu} \bigl(\alpha_1(A),\ldots,\alpha_{g}(A)\bigr).$$
It follows that computing 
\begin{equation} \label{eq-a}
a_{\mu}(\A_g;q) := \sum_{[A] \in [\A_g(\FF_q)]} \frac{a_{\mu}(A)}{|\Aut_{\FF_q}(A)|} 
\end{equation}
for each $\mu$ such that $|\mu|\leq |\lambda|$ gives a way to compute \eqref{eq-traceFq}. 

For any substack $X \subset \A_g$ we define $a_{\mu}(X;q)$ in the corresponding way. It is used repeatedly below that for any $[A] \in [\A_g(\FF_q)]$ it is enough to know $a_i(A)$ for all $1 \leq i \leq g$ to compute $a_{\mu}(A)$ for any partition $\mu$.

We now reformulate \eqref{eq-a} for genus~$3$ in terms of curves using the Torelli morphism, $t_g\colon \M_g \to \A_g.$ This morphism between stacks is of degree $1$ if $g \leq 2$ (where the genus~$1$ case should be interpreted as the isomorphism $t_1\colon \M_{1,1} \to \A_1$) and of degree $2$ if $g \geq 3$, ramified along the hyperelliptic locus $\Hh_g \subset\M_g$. 

Let $[\M_g(\FF_q)]$ denote the set of $\FF_q$-isomorphism classes of smooth curves of genus~$g$ defined over $\FF_q$. It is essential for us that if $[C] \in [\M_g(\FF_q)]$, then $H^1(C\otimes \FFb_q,\QQl)$ and $H^1(t_g(C)\otimes \FFb_q,\QQl)$ are isomorphic as $\Galp{\FF}{q}$-modules. We define $a_{\mu}(C)$ in the same way as for an abelian variety and we note that by the Lefschetz trace formula, $|C(\FF_{q^i})|=1+q^i-a_i(C)$. 

Turning to genus~$3$ and putting $\M_3^0:=\M_3\setminus \Hh_3$ and $\A_{1,1,1}:=\A_1^3/\s_3$, we have the following stratification: 
\begin{equation} \label{eq-decomp}
\A_3 = t_3(\M_3^0) \; \sqcup \; t_3(\Hh_3) \; \sqcup \; (t_2(\M_2) \times \A_1) \;\sqcup \;\A_{1,1,1}.   
\end{equation}
Therefore $a_{\mu}(\A_3;q)=a_{\mu}(t_3(\M_3^0);q)+a_{\mu}(t_3(\Hh_3);q)+a_{\mu}(t_2(\M_2) \times \A_1;q)+a_{\mu}(\A_{1,1,1};q)$, and we now turn to the computation of the different terms in this sum. 

\subsection{Non-hyperelliptic curves of genus~$3$}
For any $[C] \in [\M_3^0(\FF_q)]$, there are precisely two elements of $[\A_3(\FF_q)]$ whose representatives are isomorphic over $\FFb_q$ to the Jacobian $J(C)$, namely, the Jacobian $J(C)$ itself and its ``twist'' $J(C)^{-1}$, see for instance \cite{Serre-appendix}, Appendix.
On the other hand, the automorphism group of $J(C)$ includes the element $-1$, which doesn't come from an automorphism of $C$. It directly follows from this that $a_{\mu}(t_3(\M_3^0);q)=a_{\mu}(\M_3^0;q)$ if $|\mu|$ is even. If $|\mu|$ is odd, then $a_{\mu}(t_3(\M_3^0);q)=0$, but this does not necessarily hold for $a_{\mu}(\M_3^0;q)$. 

For any non-hyperelliptic curve $C$ of genus~$3$, the canonical linear system gives an isomorphism of $C$ with a plane quartic curve. Conversely, any non-singular plane quartic curve is non-hyperelliptic of genus~$3$. Any plane quartic curve can be given by a homogeneous degree $4$ polynomial in three variables. Identify this space of polynomials with $\PP^{14}$ and let $Q \subset \PP^{14}$ denote the subset that gives rise to non-singular quartics. All isomorphisms between plane quartic curves are induced by $\PGL_3$, the automorphism group of the plane. It follows that
\begin{equation*} 
a_{\mu}(\M_3^0;q)=\frac{1}{|\PGL_3(\FF_q)|}\sum_{C \in Q(\FF_q)} a_{\mu}(C).
\end{equation*}
For all $q \leq 25$ and $C \in Q(\FF_q)$, we have computed $|C(\FF_q)|$, $|C(\FF_{q^2})|$ and $|C(\FF_{q^3})|$. From this information, we can then compute $a_{\mu}(\M_3^0;q)$ for any partition $\mu$.

\subsection{Hyperelliptic curves} 
The Torelli morphism gives a bijection between $[\Hh_g(\FF_q)]$ and $[t_g(\Hh_g)(\FF_q)]$, and $\Aut_{\FF_q}(C)=\Aut_{\FF_q}(t_g(C))$ for every $[C] \in [\Hh_g(\FF_q)]$. It follows that $a_{\mu}(t_g(\Hh_g);q)=a_{\mu}(\Hh_g;q)$ for all partitions $\mu$.

A hyperelliptic curve of genus~$g \geq 2$ comes with a canonical separable degree two morphism to $\PP^1$. Using this morphism, we can describe any hyperelliptic curve as a separable degree two extension of $k(\PP^1) \cong k(x)$, where we have chosen a coordinate $x$. If we assume that the characteristic is not two, then the degree two extension can be written in the form $y^2=f(x)$, where $f$ is a square-free polynomial of degree $2g+1$ or $2g+2$. Denote the set of such polynomials by $P_g$. The isomorphisms between curves written in this form are generated by $\PGL_2=\Aut(\PP^1)$ and scalar multiplications of $y$. If $C_f$ is the hyperelliptic curve corresponding to $f \in P_g$, then 
$$a_i(C_f)=-\sum_{x \in \PP^1(\FF_{q^i})}\chi_{2,i}(f(x)),$$
where $\chi_{2,i}$ is the quadratic character and $f(\infty)$ is equal to the $(2g+2)$nd coefficient of $f$. So, for characteristics not equal to two:
\begin{equation*} 
a_{\mu}(\Hh_g;q)=\frac{1}{|\GL_2(\FF_q)|}\sum_{f \in P_g(\FF_q)} a_{\mu}(C_f).  
\end{equation*}
For a corresponding expression in characteristic two, see for instance the description in \cite[Section 8]{B2}. For all $q \leq 17$ and $f \in P_3(\FF_q)$, we have computed $|C_f(\FF_q)|$, $|C_f(\FF_{q^2})|$ and $|C_f(\FF_{q^3})|$. From this information, we can then compute $a_{\mu}(\Hh_3;q)$ for any partition $\mu$ and $q \leq 17$. 

An elliptic curve comes with a marked point and it has a canonical separable degree two morphism to $\PP^1$ such that the marked point is a ramification point over infinity. Thus, the elliptic curves can also be written in the form $y^2=f(x)$, where $f$ belongs to the set $P'_1$ of square-free polynomials of degree~$3$, and the group $G$ of isomorphisms is generated by scalar multiplications of $y$ and morphisms induced from elements of $\PGL_2$ that keep infinity fixed.

\subsection{The decomposable abelian threefolds}
Consider first the abelian threefolds that are isomorphic to a product of an elliptic curve and a Jacobian of a genus~$2$ curve. The first cohomology group of such an abelian variety is equal to the direct sum of the first cohomology groups of the two curves. Moreover, the automorphisms of the abelian variety come only from the automorphisms of the curves. We therefore find that $a_{\mu}(t_2(\M_2) \times \A_1 ;q)$ equals
$$\frac{1}{|\GL_2(\FF_q)|}\frac{1}{|G(\FF_q)|} \sum_{f \in P_2(\FF_q)} \sum_{h \in P'_1(\FF_q)} \prod_{i=1}^{\mu_1} \Bigl( a_{\hat \mu_i}(C_f)+ a_{\hat \mu_i}(C_h) \Bigr).$$
From the computation, for all $q \leq 17$, of $|C_f(\FF_q)|$ and $|C_f(\FF_{q^2})|$ for all $f \in P_2(\FF_q)$ and $|C_h(\FF_q)|$ for all $h \in P'_1(\FF_q)$, we then compute $a_{\mu}\bigl(t_2(\M_2) \times \A_1 ;q \bigr)$ for any partition $\mu$. 

Let us now consider an abelian threefold that is a product of three elliptic curves. Again, the first cohomology group of such an abelian threefold is the direct sum of the first cohomology groups of the curves. There are three possibilities for the abelian threefold to be defined over $\FF_q$. Either all three elliptic curves are defined over $\FF_q$; or one of them is defined over $\FF_q$ and the two others are defined over $\FF_{q^2}$, where one is sent to the other by Frobenius $F_q$; or finally, all are defined over $\FF_{q^3}$ and Frobenius $F_q$ permutes these three curves cyclically. If we are in the first case and any of the curves are equal, then there are extra automorphisms coming from the possible permutations. In the second case, $F_q$ becomes a new non-trivial automorphism, and in the third case we have both $F_q$ and $F_q^2$. These observations leave us with the following formula for $a_{\mu}\bigl( \A_{1,1,1} ;q \bigr)$:
\begin{multline*}
\frac{1}{6}\, \frac{1}{|G(\FF_q)|^3} \sum_{f_1,f_2,f_3 \in P'_1(\FF_q)} \prod_{i=1}^{\mu_1} \bigl( a_{\hat \mu_i}(C_{f_1})+a_{\hat \mu_i}(C_{f_2})+a_{\hat \mu_i}(C_{f_3}) \bigr)+{}\\  
+\frac{1}{2}\, \frac{1}{|G(\FF_q)|}\frac{1}{|G(\FF_{q^2})|}\sum_{f \in P'_1(\FF_q)}\sum_{h \in P'_1(\FF_{q^2})} \prod_{i=1}^{\mu_1} \bigl( a_{\hat \mu_i}(C_f)+a_{\hat \mu_i}(C_h) \bigr)+{}\\
+\frac{1}{3}\, \frac{1}{|G(\FF_{q^3})|}\sum_{f \in P'_1(\FF_{q^3})} \prod_{i=1}^{\mu_1} \bigl( a_{\hat \mu_i}(C_f) \bigr).
\end{multline*}
For all $q \leq 17$, $1 \leq j \leq 3$ and $f_j \in P'_1(\FF_{q^j})$, we have computed $|C_{f_j}(\FF_{q^j})|$. From this information, we can then compute $a_{\mu}(\A_{1,1,1};q)$ for any partition $\mu$ and $q \leq 17$. 

\subsection{Closed formulas for numbers of points}
In the articles~\cite{B1} and \cite{B2}, the first author counted points over \emph{any} finite field for the spaces $\M^{0}_{3,n}$ and $\Hh_{g,n}\,$, for $g \geq 2$ and $n \leq 7$. Using this information, $e_c(\M^{0}_3 \otimes \overline{\QQ},\VVl)$ and $e_c(\Hh_g \otimes \overline{\QQ},\VVl)$, for all $g \geq 2$, could be determined as elements of $K_0(\GalC)$ for all $|\lambda| \leq 7$. 
This information can be pieced together (with genus~$1$ information from Theorem~\ref{thm-g1}) using the decomposition \eqref{eq-decomp}, compare \cite{JBvdG}, to determine  
$e_c(\A_3 \otimes \FFb_q,\VVl)$ for $|\lambda| \leq 6$ as an element of $K_0(\GalC)$. 
The results are polynomials in $\QQl(-1)$.

\begin{theorem} Let $\qq$ denote the class of $\QQl(-1)$ in $K_0(\GalC)$. The following holds: 

\bigskip
\vbox{
\centerline{\def\quad{\hskip 0.3em\relax}
\vbox{\offinterlineskip
\hrule
\halign{&\vrule#& \quad \hfil#\hfil \strut \quad  \cr
height2pt&\omit&&\omit&&\omit&&\omit& \cr 
& $\lambda$ && $e_c(\A_{3}\otimes \overline{\QQ},\VVl)$ && $\lambda$ && $e_c(\A_{3}\otimes \overline{\QQ},\VVl)$  & \cr
height2pt&\omit&&\omit&&\omit&&\omit& \cr 
\noalign{\hrule}
height2pt&\omit&&\omit&&\omit&&\omit& \cr 
&$(0,0,0)$ && $\qq^6+\qq^5+\qq^4+\qq^3+\mathbf{1} $ && 
$(2,0,0)$ && $-\qq^3-\qq^2 $  &\cr
height2pt&\omit&&\omit&&\omit&&\omit& \cr 
&$(1,1,0)$ && $-\qq $ && 
$(4,0,0)$ && $-\qq^3-\qq^2 $  &\cr
height2pt&\omit&&\omit&&\omit&&\omit& \cr 
& $(3,1,0)$ && $0 $ && 
$(2,2,0)$  && $0 $  &\cr
height2pt&\omit&&\omit&&\omit&&\omit& \cr 
& $(2,1,1)$ && $\mathbf{1} $ && 
$(6,0,0)$ && $-2\qq^3-\qq^2 $  &\cr
height2pt&\omit&&\omit&&\omit&&\omit& \cr 
&$(5,1,0)$ && $-\qq^4 $ && 
$(4,2,0)$ && $-\qq^5+\qq $  &\cr
height2pt&\omit&&\omit&&\omit&&\omit& \cr 
&$(4,1,1)$ && $\mathbf{1} $ && 
$(3,3,0)$ && $\qq^7-\qq $  &\cr
height2pt&\omit&&\omit&&\omit&&\omit& \cr 
&$(3,2,1)$ && $0 $ && 
$(2,2,2)$ && $\mathbf{1} $  &\cr
height2pt&\omit&&\omit&&\omit&&\omit& \cr 
} \hrule}
}}
\end{theorem}

\begin{remark} If we assume it to be known that $s_{n(\lambda)}=0$ for all $\lambda$ with $|\lambda| \leq 6$, then this theorem proves that Conjecture~\ref{conj-g3} is true for these local systems. 
\end{remark}

%%%%%%%%%%%%%%%%%%%%%%%%%%%%%%%%%%%%
\section{Characteristic polynomials}  \label{sec-g3-charpol} 
In this section, we will use the conjectures stated in Section~\ref{sec-g3} together with the information coming from our counts of points over finite fields to (conjecturally) determine the local spinor $L$-factor at $p=2$ of certain Siegel modular forms of degree $3$.

Let $f\in \Sigma^{\gen}_{a-b,b-c,c+4}$ be an eigenform and put $m=(a+b+c+6)/2$. The local spinor $L$-factor $L_p(s,f)$ equals $Q_p(p^{-s},f)^{-1}$, where $Q_p(X,f)$ is the characteristic polynomial of the image of the conjugacy class $s_p(f)$ in $\mathrm{GSpin}_{7}(\CC)$ under the spin representation. It is known that $Q_p(X,f)$ satisfies the following duality:
$$Q_p(X,f)=(p^{m}X)^8Q_p(p^{-2m}X^{-1},f).$$

Let $\Lambda$ be the set of $\lambda$ for which Conjecture~\ref{conj-g3-int} together with Conjecture~\ref{conj-g3-lifts} indicate that $s^{\gen}_{\lambda}=1$, that is, that the space of generic Siegel modular forms is one-dimensional (with $|\lambda| \leq 60$). This set consists of $123$ elements. For $\lambda$ in $\Lambda$, a single eigenform $f_{\lambda}$ should generate the space of generic Siegel modular forms.

We thus expect $\Sigma^{\gen}[n(\lambda)]$ to be of rank $8$ for $\lambda$ in $\Lambda$. Assuming this, we can use Conjecture~\ref{conj-g3} and Conjecture~\ref{conj-g3-nonexh} together with the computations of $\Trace(F_{2^i},e_c(\A_3 \otimes \FFb_{2^i},\VVl))$ for $1 \leq i \leq 4$ to compute the first five coefficients of a polynomial $F_{2,\lambda}(X)$ of degree $8$, which should equal the characteristic polynomial $C_2(X,\Sigma^{\gen}[n(\lambda)])$ of $F_2$. Recall from Section~\ref{sec-motives} that $C_2(X,\Sigma^{\gen}[n(\lambda)])$ should in turn equal $Q_2(X,f_{\lambda})$. Thus, using the duality above, we have a conjectural way of finding the local spinor $L$-factor at $p=2$ for $123$ Siegel modular forms of degree $3$.

We will now make two ``checks'' of the polynomials $F_{2,\lambda}(X)$ that we have found. We view these as substantial consistency checks for our conjectures.

First, for all $\lambda$ in $\Lambda$, we find that the Ramanujan conjecture is fulfilled for $F_{2,\lambda}(X)$, that is, the roots of $F_{2,\lambda}(X)$ have absolute value $2^{-m}$.

If we normalize $Q_p(X,f)$ by 
$$\hat Q_p(X,f):= Q_p(p^{-m}X,f),$$
we will get a polynomial of the form 
\begin{equation} \label{eq-charpol}
X^8-AX^7+BX^6-CX^5+DX^4-CX^3+BX^2-AX+1,
\end{equation}
for which (see \cite[Prop.~2.2.2]{Gross-Savin})
\begin{equation} \label{eq-cond-Spin7}
A^2(D+2B+1)=C^2+2AC+A^4.
\end{equation}
If we fix $\lambda$ in $\Lambda$ and put 
$$\hat F_{2,\lambda}(X):=F_{2,\lambda}(2^{-m}X),$$
we get a polynomial of the form \eqref{eq-charpol}. The second check is that for all $\lambda$ in $\Lambda$ the condition \eqref{eq-cond-Spin7} holds for $\hat F_{2,\lambda}(X)$.

\subsection{Motives for $\Gtwo$} 
In the article \cite{Serre-galmot}, Serre asked if there are motives with motivic Galois group of type $\Gtwo$. Gross and Savin suggested in \cite{Gross-Savin} that one should search for such motives inside the cohomology of a moduli space of principally polarized abelian varieties of dimension~$3$ with an Iwahori level structure at a finite number of primes. They were able to construct automorphic forms of level $2$ and of level $5$ for $\mathrm{PGSp}_6$ that both are lifts from automorphic forms for $\Gtwo$ (\cite[Prop.~4.5.8]{Gross-Savin}). 

The dual group of an anisotropic form of $\Gtwo$ equals $\Gtwo(\CC)$ and we have an inclusion of dual groups $\Gtwo(\CC) \hookrightarrow \mathrm{Spin}_{7}(\CC)$, which realizes $\Gtwo(\CC)$ as the stabilizer of a non-isotropic vector in the spin representation, see \cite{Gross-Savin}. The Langlands program then predicts that there should be a lift of automorphic forms from $\Gtwo$ to $\mathrm{GSp}_6$. For a Siegel modular eigenform $f$ for $\mathrm{GSp}_6$ that is a lift from $\Gtwo$, the normalized spinor $L$-function $\hat Q_p(X,f)$, which is of the form \eqref{eq-charpol}, will fulfil the relation \eqref{eq-cond-Spin7} and the relation
\begin{equation} \label{eq-cond-G2}
2A-2B+2C-D-2=0.
\end{equation}
These relations have the consequence that $\hat Q_p(X,f)$ has $(X-1)^2$ as a factor. 

We expect the motive corresponding to a Siegel modular eigenform $f$ that is a lift from $\Gtwo$ to have rank $8$, and it should decompose into an irreducible piece of dimension $7$ and a Tate class. The Tate class will leave a footprint on $\hat Q_p(X,f)$ as one of the factors $(X-1)$ mentioned above. In \cite{Gross-Savin}, it is shown that if there are lifts from $\Gtwo$ to $\mathrm{GSp}_6$, then they will appear for local systems of the form $(b+c,b,c)$. 

For all $\lambda$ in $\Lambda$, we can compute $\hat F_{2,\lambda}(X)$ as described above, and we found $\hat F_{2,\lambda}(X)$ to be irreducible except in the following $12$ cases:
\begin{multline*}
 (9, 6, 3), (10, 8, 2), (10, 6, 4), (11, 10, 1), (11, 8, 3), (11, 7, 4), (12, 12, 0),\\ (12, 9, 3), (12, 6, 6), (14, 14, 0), (12, 11, 5), (15, 15, 0).
\end{multline*}
In all these $12$ cases, conditions \eqref{eq-cond-Spin7} and \eqref{eq-cond-G2} were found to hold. We see that all these local systems are of the form $(b+c,b,c)$, except for $\lambda=(12, 11, 5)$. For the single case, the behaviour at $p=2$ may be an anomaly. In the other $11$ cases, we expect that there is a motive of rank $7$ with motivic Galois group of type $\Gtwo$ appearing in $H^6_c(\A_3,\VVl)$.

\begin{example} For $\lambda=(9,6,3)$, we use the method above to compute
\begin{multline*}
F_{2,\lambda}(X)=(1-2^{12}\, X)^2 \, (1+7112 \, X+34431488 \, X^2+\\+176085008384 \, X^3+2^{24} \cdot 34431488 \, X^4+2^{48} \cdot 7112 \, X^5+2^{72} \, X^6),
\end{multline*}
which we expect to be the characteristic polynomial at $p=2$ of a modular form lifted from $\Gtwo$.
\end{example}

%%%%%%%%%%%%%%%%%%%%%
\section{Congruences} \label{sec-congruences}
Motivated by his work on the Eisenstein cohomology of a local system on ${\A}_2$, Harder formulated in \cite{Ha,Ha1} a conjectural congruence between the Hecke eigenvalues of an elliptic cusp form and those of a degree $2$ Siegel modular form modulo a `large' prime that divides a critical value of the $L$-function of the elliptic modular form. In this section, we will present this conjecture and several similar conjectural congruences for Siegel modular forms of degree~$2$ and $3$. These conjectures were formulated in close collaboration with Harder. 

We view the conjectures in this section as the beginnings of a general theory
of congruences, and every concrete example of a congruence as an additional
confirmation of the conjectures in sections 6 and 7, for which we have much
more evidence than for the conjectures below. As the case of Harder's
conjecture shows, special care will be required 
when the local system given by $(a,b)$ or $(a,b,c)$ is not regular.

\begin{definition} 
For a Siegel modular eigenform $f$ of degree $g$, a prime number~$p$, and
an integer $r \geq 1$, we use the Satake parameters to define
$$
\lambda_{p^r}(f)=
\sum_{s=0}^g\, 
\sum_{1\leq i_1 < i_2 < \ldots 
<i_s\leq g} \bigl(\alpha_{p,0}(f)\alpha_{p,i_1}(f) \cdots \alpha_{p,i_s}(f)\bigr)^r.
$$
\end{definition}
\noindent For $r=1$, this definition coincides with our previous one: the eigenvalue of $f$ under the action of the Hecke operator $T(p)$. 
\begin{remark} Note that this notation is non-standard.
\end{remark}

\begin{definition} For a Siegel modular eigenform $f$, we have a finite field extension, namely $\QQ_{f}:=\QQ(\lambda_p(f): p \text{ prime})$. If $f_1,\ldots,f_n$ are eigenforms, then $\QQ_{f_1,\ldots,f_n}$ will denote the compositum of the fields $\QQ_{f_1},\ldots,\QQ_{f_n}$.
\end{definition}

For an eigenform $f \in S_{k}$, we put $L_{\infty}(f,s):=\Gamma(s)/(2\pi)^s$ and $\Lambda(f,s):=L_{\infty}(f,s) \, L(f,s)$. The function $\Lambda(f,s)$ has a holomorphic continuation to the whole complex plane and it fulfils the functional equation $\Lambda(f,s)=(-1)^{k/2}\Lambda(f,k-s)$. The numbers $\Lambda(f,r)$ for $k/2 \leq r \leq k-1$ are called \emph{critical values} of $f$. By a theorem of Manin and Vishik, see \cite{Haberland}, there exist two real numbers $\omega_{+}(f)$ and $\omega_{-}(f)$, called \emph{periods}, such that $\Lambda'(f,r):=\Lambda(f,r)/\omega_{\pm}(f)$ (where we take $\omega_{+}$ if $r$ is even, and $\omega_{-}$ if $r$ is odd) lies in the field $\QQ_{f}$. We will say that a prime $\ell$ in $\QQ_{f}$ divides $\Lambda'(f,r)$ if it divides the numerator of $\Lambda'(f,r)$. Moreover, if $\ell$ in $\QQ_{f}$ lies above the prime $p$ in $\QQ$, then $\ell$ will be called ordinary (for $f$) if $\lambda_p(f) \equiv_{\ell} 0$ does not hold. 

\begin{conjecture}[Harder \cite{Ha,vdG1}] \label{conj-cong-harder} Assume $a>b$. Take an eigenform $f \in S_{a+b+4}$. If for an ordinary prime $\ell$ in $\QQ_{f}$ and $s\geq 1$ the number $\ell^s$ divides the critical value $\Lambda'(f,a+3)$, then there is an eigenform $F \in S_{a-b,b+3}$ such that
\begin{equation*} 
\lambda_{q}(F) \equiv_{\ell^s} \lambda_{q}(f)+q^{a+2}+q^{b+1} 
\end{equation*}
holds in $\QQ_{f,F}$ for all prime powers $q$.
\end{conjecture}

\begin{remark} The original conjecture supposed $\ell$ to be `large' (a slightly imprecise notion, but $\ell$ greater than the weight of the form should suffice). Harder suggested to replace it by `ordinary.' The conjectural congruence has been checked numerically in many  instances, see \cite[p.~237--240]{vdG1}. 

The congruences in Harder's conjecture should come from denominators of
Eisenstein classes in the Betti cohomology. In the case $a=b$ even, there is no
Eisenstein class (see \cite{Ha2}) and hence no conjectured congruence either.
In the case $a=b$ odd, the conjecture is trivially true due to the presence
of Saito-Kurokawa lifts. 
The statement with $S_{a-b,b+3}$ replaced by
$\Sigma^{\gen}_{a-b,b+3}$ is nontrivial.
It has been proved under certain conditions, see \cite[Section 6]{Katsurada-SK}
and \cite{DIK}.
\end{remark}

For any elliptic modular eigenforms $f_1,\ldots,f_m$, where $f_i \in S_{k_i}\,$, we define the following $L$-function through an Euler product:
$$L \bigl(\bigotimes_{i=1}^m \mathrm{Sym}^{r_i}(f_i) ,s 
\bigr):=\prod_{p} \Bigl( \prod_{i=1}^m \prod_{j=0}^{r_i} 
\bigl(1-\alpha_{p,0}(f_i)^{r_i} \, \alpha_{p,1}(f_i)^{j} \, p^{-s} 
\bigr) \Bigr)^{-1}.$$

Below, several different instances, here denoted $L (\cdot ,s)$, of this $L$-function will appear and we will naively assume that we have corresponding properties to what we saw above for $L(f,s)$ (for more details, see \cite{Deligne-L}). That is, that we can define a factor at infinity $L_{\infty}(\cdot,s)$, such that $\Lambda(\cdot,s):=L_{\infty}(\cdot,s) \, L(\cdot,s)$ 
has a meromorphic continuation to the complex plane and fulfils a functional equation.
Moreover, that there is a set of integers $r$ called critical values and two real numbers $\omega_{+}(f)$ and $\omega_{-}(f)$, such that for each critical value $r$, the quotient $\Lambda'(\cdot,r):=\Lambda(\cdot,r)/\omega_{\pm}(f)$ lies in the field $\QQ_{f}$. We will then, in the same way as above, talk about primes dividing $\Lambda'(\cdot,r)$. In the two cases appearing in the section for genus~$2$, these properties are known to hold; see \cite{Zagier-mod} for the case $L(\Sym^2 f,s)$. The conjectures we will formulate will be slightly imprecise in the sense that we will not specify what is meant by a \emph{large} prime. One may suspect that $\ell$ being ordinary for the $g=1$ forms involved should suffice.

Our general philosophy can somewhat vaguely be formulated as follows. Congruences appear between elements of different subspaces of $H^i_c(\A_g,\VVl)$ invariant under the Hecke algebra, due to integrality issues (compare \cite[p.~255]{Ha}), and their appearance is governed by vanishings of critical values of $L$-functions modulo primes. We are interested in congruences involving generic Siegel modular forms of genus $g$ and thus we consider the middle cohomology group. 

In this picture, the following congruence in genus $1$ is connected to the class in $H^1_{\Eis}(\A_1,\VV_{a})$ that comes from an Eisenstein series whose Hecke eigenvalue for $T(p)$ equals $1+p^{a+1}$. If a prime $\ell$ divides the numerator of $\zeta(-a-1)$, then there is an eigenform $F \in S_{a+2}$ such that
$$\lambda_q(F) \equiv_{\ell} 1+q^{a+1}$$
holds in $\QQ_{F}$ for all prime powers $q$. 

In a similar way, Harder's congruence should be connected to the classes $s_{a+b+4}L^{b+1}$, found in $H^3_{\Eis}(\A_2,\VV_{a,b})$, if $a \neq b$, and to the classes coming from Saito-Kurokawa lifts, if $a=b$ is odd (cf.~\cite{Ha2}).

\subsection{Congruences in genus~2}
First we put 
\begin{equation*}
e_{2,\rm SK}(a,b):=-S[2a+4]-s_{2a+4}(L^{a+1}+L^{a+2}) \quad \text{if $a=b$ is odd} 
\end{equation*}
and then 
$$N_q(a,b):=-\Trace \bigl(F_q,e_c({\A}_2,\VV_{a,b})\bigr)+\Trace \bigl(F_q,e_{2,\rm extr}(a,b) \bigr)+\Trace\bigl(F_q,e_{2,\rm SK}(a,b)\bigr).$$

\subsubsection{Kurokawa-Mizumoto congruence}
The following conjecture was formulated together with Harder (it should be compared with \cite[Prop. 4.4]{Dummigan-sym-II}) and is connected to the Eisenstein series which give rise to the classes $S[a+3]$ found in $H^3_{\Eis}(\A_2,\VV_{a,b})$. It generalizes, to the vector-valued case, a congruence found by Kuroka\-wa, which was proved (also in the case of higher genera) by Katsurada and Mizumoto, see \cite{Katsurada-Mizumoto} and \cite{Mizumoto}. The Kuroka\-wa-Mizumoto congruence is the case $a=b$ in our conjecture. Several instances of this generalized congruence have already been proved for vector-valued Siegel modular forms, see \cite{Satoh} and \cite{Dummigan-sym-II}.

\begin{conjecture} \label{conj-Mg2cong} Take an eigenform $f \in S_{a+3}$. If for a ``large'' prime $\ell$ in $\QQ_{f}$ and $s\geq 1$ the number $\ell^s$ divides the critical value $\Lambda'(\Sym^2(f),a+b+4)$, then there is an eigenform $F \in \Sigma^{\gen}_{a-b,b+3}$ such that
\begin{equation} \label{eq-Mg2cong} 
\lambda_{q}(F) \equiv_{\ell^s} \lambda_{q}(f) \, \bigl(q^{b+1}+1 \bigr)
\end{equation}
holds in $\QQ_{f,F}$ for all prime powers $q$.
\end{conjecture}

Let $f \in S_k$ and $k \leq 22$ be such that a ``large'' prime $\ell$ to the power $s$ divides $\Lambda'(\Sym^2(f),m)$ for $k+1 \leq m \leq 2k-2$, see for instance \cite[Table~1]{Dummigan-sym}. We have then checked that $s^{\gen}_{a-b,b+3} \geq 1$ for the corresponding pairs $(a,b)$. There are seven cases for which $s^{\gen}_{a-b,b+3}=1$ and in all these, $N_q(a,b)$ is congruent, modulo $\ell^s$, to the right hand side of \eqref{eq-Mg2cong} for all $q \leq 37$. 

\subsubsection{Yoshida-type congruence} \label{sec-yoshida}
This congruence should be connected to the endoscopic contribution $-s_{a+b+4} S[a-b+2] L^{b+1}$, and it was found following a suggestion by Zagier. 

\begin{conjecture} \label{conj-Yg2cong} Take any eigenforms $f \in S_{a+b+4}$ and $g \in S_{a-b+2}$. If for a ``large'' prime $\ell$ in $\QQ_{f}$ and $s\geq 1$ the number $\ell^s$ divides the critical value $\Lambda'(f \otimes g,a+3)$, then there is an eigenform $F \in S_{a-b,b+3}$ such that
\begin{equation} \label{eq-Yg2cong} 
\lambda_{q}(F) \equiv_{\ell^s} \lambda_{q}(f)+q^{b+1} \lambda_{q}(g)
\end{equation}
holds in $\QQ_{f,g,F}$ for all prime powers $q$.
\end{conjecture}

That $a+3$ was the correct critical value of the $L$-function was suggested by Dummigan. He computed (algebraically) the critical value $\Lambda'(f \otimes g,a+3)$ for eigenforms $f \in S_{a+b+4}$ and $g \in S_{a-b+2}$, when $(a,b)=(18,8)$, $(20,4)$, $(21,3)$ or $(20,0)$. The only large primes dividing the norm of these critical values were $\ell=263$ for $(18,8)$, $\ell=223$ for $(20,4)$, and $\ell=2747$ for $(21,3)$. In all these cases, $s_{a-b,b+3}=1$ and $N_q(a,b)$ is congruent, modulo $\ell$, to the norm of the right hand side of \eqref{eq-Yg2cong} for all $q \leq 37$. We have that $s_{20,3}=0$, but there is at the same time no large prime dividing the critical value for $(a,b)=(20,0)$.

\subsection{Congruences in genus~3}
First we put 
%\begin{multline*}
$$ \begin{aligned}
e_{3,\rm ne}(a,b,c):=S[b+3]\bigl(S[a+c+5]+L^{c+1}S[a-c+3]\bigr)&\\
{}+S[a+4]\bigl(S[2b+4]+s_{2b+4}(L^{b+1}+L^{b+2}) \bigr) \qquad &\mbox{if}\,\, {b=c} \\
{}+S[c+2]\bigl(S[2a+6]+s_{2a+6}(L^{a+2}+L^{a+3}) \bigr) \qquad &\mbox{if}\,\, {a=b} 
%+\begin{cases} S[a+4]\bigl(S[2b+4]+s_{2b+4}(L^{b+1}+L^{b+2}) \bigr)
%&\text{if $\, b=c$} \\
%S[c+2]\bigl(S[2a+6]+s_{2a+6}(L^{a+2}+L^{a+3}) \bigr)
%&\text{if $\, a=b$} \end{cases}
%\end{multline*}
\end{aligned} $$
and then 
\begin{multline*}
N_q(a,b,c):=\Trace \bigl(F_q,e_c({\A}_3,\VV_{a,b,c})\bigr) - \Trace \bigl(F_q,e_{3,\rm extr}(a,b,c) \bigr)\\-\Trace \bigl(F_q,e_{3,\rm ne}(a,b,c) \bigr).
  \end{multline*}

\subsubsection{Congruences of Eisenstein type}
The following two conjectures were formulated together with Harder and should be connected to contributions to the Eisenstein cohomology of the form $s_{b+c+4}L^{c+1}S[a+4]$ (respectively $s_{a+b+6}L^{b+2}S[c+2]$) and to non-exhaustive lifts of the form (ii) and (iii) in Conjecture~\ref{conj-g3-lifts}. 

\begin{conjecture} \label{conj-firstg3cong1} Take any eigenforms $f \in S_{a+4}$ and $g \in S_{b+c+4}$. If for a ``large'' prime $\ell$ in $\QQ_{f,g}$ and $s\geq 1$ the number  $\ell^s$ divides the critical value $\Lambda'(\Sym^2(f) \otimes g,a+b+6)$, then there is an eigenform $F \in \Sigma^{\gen}_{a-b,b-c,c+4}$ such that
\begin{equation} \label{eq-firstg3cong1} 
\lambda_{q}(F) \equiv_{\ell^s} \lambda_{q}(f) \bigl(\lambda_{q}(g)+q^{b+2}+q^{c+1} \bigr) 
\end{equation}
holds in $\QQ_{f,g,F}$ for all prime powers $q$.
\end{conjecture} 

\begin{remark} This congruence was formulated in the scalar-valued case (generalized to any Miyawaki-Ikeda lift) by Katsurada, \cite{Katsurada-letter}. 
\end{remark}

Mellit computed (using numerical approximations) a list of critical values $\Lambda'(\Sym^2(f) \otimes g,\nu)$, for all eigenforms $f \in S_{k_1}$ and $g \in S_{k_2}$ with $k_1,k_2 \leq 20$. There are $17$ cases in this list, presented in Table~\ref{tab-con}, for which a large prime power $\ell^s$ divides a critical value and (conjecturally) $s^{\gen}_{a-b,b-c,c+4}=1$, for the corresponding tuple $(a,b,c)$. In all these cases, $N_q(a,b,c)$ is congruent, modulo $\ell^s$, to the right hand side of \eqref{eq-firstg3cong1} for all $q \leq 17$. 

The critical values in Table~\ref{tab-con} have also been checked algebraically by Katsurada, and recently, Poor and Yuen have proved that the congruence holds in the case $(a,b,c)=(12,12,12)$, \cite{Katsurada-letter}.

\begin{table}[ht]\caption{Congruences of the form \eqref{eq-firstg3cong1}.}  \label{tab-con} 
\vbox{
\centerline{\def\quad{\hskip 0.3em\relax}
\vbox{\offinterlineskip
\hrule
\halign{&\vrule#& \quad \hfil#\hfil \strut \quad \cr
height2pt&\omit&&\omit&&\omit&&\omit&\cr
& $(a,b,c)$ && $\ell^s$ && $(a,b,c)$ && $\ell^s$ &\cr
\noalign{\hrule}
height2pt&\omit&&\omit&&\omit&&\omit&\cr
  & $(12,6,2)$ && $101$  && $(14,7,1)$ &&  $17^2$ & \cr
  & $(16,8,0)$ && $43$  && $(16,7,1)$ &&  $263$ & \cr
  & $(16,5,3)$ && $127$ && $(16,4,4)$ &&  $29$ & \cr
  & $(12,12,0)$ && $37$ && $(12,9,3)$ &&  $137$ & \cr
  & $(12,6,6)$ &&  $229$ && $(14,11,1)$ &&  $37$ & \cr
  & $(14,7,5)$ &&  $71$ && $(12,8,6)$ &&  $73$ & \cr
  & $(14,14,0)$ &&  $59$ && $(12,8,8)$ &&  $61$ & \cr
  & $(12,12,12)$ &&  $107$ && $(14,13,11)$ &&  $41$ & \cr
  & $(16,16,16)$ &&  $157$ &&  &&   & \cr
} \hrule}
}}
\end{table}

\begin{conjecture} \label{conj-firstg3cong2} Take any eigenforms $f \in S_{c+2}$ and $g \in S_{a+b+6}$. If for a ``large'' prime $\ell$ in $\QQ_{f,g}$ and $s\geq 1$ the number $\ell^s$ divides the critical value $\Lambda'(\Sym^2(f) \otimes g, a+c+5)$, then there is an eigenform $F \in \Sigma^{\gen}_{a-b,b-c,c+4}$ such that
\begin{equation} \label{eq-firstg3cong2} 
\lambda_{q}(F) \equiv_{\ell^s} \lambda_{q}(f)\bigl(\lambda_{q}(g)+q^{a+3}+q^{b+2}\bigr) 
\end{equation}
holds in $\QQ_{f,g,F}$ for all prime powers $q$.
\end{conjecture}

We only have one example of this congruence, namely, when $(a,b,c)=(13,11,10)$. There are then eigenforms  $f \in S_{c+2}$ and $g \in S_{a+b+6}$ for which $\ell=199$ divides the norm of the critical value $\Lambda'(\Sym^2(f) \otimes g, a+c+5)$, computed (numerically) by Mellit. Conjecturally, $s^{\gen}_{a-b,b-c,c+4}=1$, and we have that $N_q(a,b,c)$ is congruent, modulo $\ell$, to the norm of the right hand side of \eqref{eq-firstg3cong2} for all $q \leq 17$.

There is also a contribution to the Eisenstein cohomology of the form $s_{a+b+6}L^{b+2}S[a-b+2]$, which should be connected to a congruence of the form
\begin{equation} \label{eq-firstg3cong3}
  \lambda_{q}(F) \equiv_{\ell^s} \bigl(\lambda_{q}(f)+q^{b+1}\lambda_{q}(g)\bigr) (1+q^{c+1})
\end{equation}
for $F \in S_{a-b,b-c,c+4}$, $f \in S_{a+b+6}$ and $g \in S_{a-b+2}$. We have three examples of $N_q(a,b,c)$ being congruent, modulo a ``large'' prime, to the right hand side of \eqref{eq-firstg3cong3} for all $q \leq 17$, and they are presented in Table~\ref{tab-con3}.

\begin{table}[ht]\caption{Congruences of the form \eqref{eq-firstg3cong3}.}  \label{tab-con3} 
\vbox{
\centerline{\def\quad{\hskip 0.3em\relax}
\vbox{\offinterlineskip
\hrule
\halign{&\vrule#& \quad \hfil#\hfil \strut \quad \cr
height2pt&\omit&&\omit&&\omit&&\omit&\cr
& $(a,b,c)$ && $\ell^s$ && $(a,b,c)$ && $\ell^s$ &\cr
\noalign{\hrule}
height2pt&\omit&&\omit&&\omit&&\omit&\cr
  & $(14,4,2)$ && $103$  && $(15,5,4)$ &&  $691$ & \cr
  & $(24,2,2)$ && $31$  &&  &&  & \cr
} \hrule}
}}
\end{table}

Finally, we have a contribution to the Eisenstein cohomology of the form $S[a-b,b+4]$, which should be connected to a congruence of the form
\begin{equation} \label{eq-firstg3cong4}
  \lambda_{q}(F) \equiv_{\ell^s} \lambda_{q}(f)(1+q^{c+1})
\end{equation}
for $F \in \Sigma^{\gen}_{a-b,b-c,c+4}$ and $f \in \Sigma^{\gen}_{a-b,b+4}$. We have found twelve examples of $N_q(a,b,c)$ being congruent, modulo a ``large'' prime, to the right hand side of \eqref{eq-firstg3cong4} for all $q \leq 17$, and they are presented in Table~\ref{tab-con4}.

\begin{table}[ht]\caption{Congruences of the form \eqref{eq-firstg3cong4}.}  \label{tab-con4} 
\vbox{
\centerline{\def\quad{\hskip 0.3em\relax}
\vbox{\offinterlineskip
\hrule
\halign{&\vrule#& \quad \hfil#\hfil \strut \quad \cr
height2pt&\omit&&\omit&&\omit&&\omit&\cr
& $(a,b,c)$ && $\ell^s$ && $(a,b,c)$ && $\ell^s$ &\cr
\noalign{\hrule}
height2pt&\omit&&\omit&&\omit&&\omit&\cr
  & $(12,6,2)$ && $149$  && $(10,6,4)$ &&  $41$ & \cr
  & $(13,5,4)$ && $601$  && $(12,8,2)$ && $59$ & \cr
  & $(12,6,4)$ && $379$  && $(20,2,2)$ && $37$ & \cr
  & $(15,5,4)$ && $29$  && $(13,7,4)$ && $23$ & \cr
  & $(13,7,6)$ && $1621$  && $(12,8,6)$ && $53$ & \cr
  & $(12,8,8)$ && $89$  && $(16,16,16)$ && $691$ & \cr
} \hrule}
}}
\end{table}

In our picture, there are two other possibilities, which are variants of the
two latter congruences, but where the contributions are found in
$e_c(\A_2,\VV_{b,c})$ instead of in $e_c(\A_2,\VV_{a+1,b+1})$. But we have
not found any examples.

\subsubsection{Congruences of endoscopic type } \label{sec-g3second}
In the following two congruences, the critical value $a+b+5$ of the $L$-function was suggested by Dummigan, on the grounds that it would (like all the congruences in this section) fit well with the Bloch-Kato conjecture. The first congruence is between lifts of the form (i) in Conjecture~\ref{conj-g3-lifts} and forms in $\Sigma^{\gen}_{a-b,b-c,c+4}$.
\begin{conjecture} \label{conj-secondg3cong1} Take any eigenforms $f \in 
S_{b+3}$, $g \in S_{a+c+5}$ and $h \in S_{a-c+3}$. If for 
a ``large'' prime $\ell$ in $\QQ_{f,g}$ and $s\geq 1$ the number $\ell^s$ divides the critical value $\Lambda'(\mathrm{Sym}^2(f) \otimes g \otimes h,a+b+5)$, then there is an eigenform $F \in \Sigma^{\gen}_{a-b,b-c,c+4}$ such that
\begin{equation} \label{eq-secondg3cong1} 
\lambda_{q}(F) \equiv_{\ell^s} \lambda_{q}(f)\bigl(\lambda_{q}(g)+q^{c+1} \lambda_{q}(h)\bigr)
\end{equation}
holds in $\QQ_{f,g,h,F}$ for all prime powers $q$.
\end{conjecture}

Mellit computed the critical value $\Lambda'(\Sym^2(f) \otimes g \otimes h,a+b+5)$ when $(a,b,c)=(12,9,3)$ and found that it was divisible by $\ell=37$. We have conjecturally that $s^{\gen}_{a-b,b-c,c+4}=1$, and $N_q(a,b,c)$ is congruent, modulo $\ell$, to the right hand side of \eqref{eq-secondg3cong1} for all $q \leq 17$. 

The following congruence should be connected to the endoscopic contribution of the form $S[a+4]\bigl(S[b+c+4]+L^{c+1}S[b-c+2 \bigr)$.
\begin{conjecture} \label{conj-secondg3cong2} Take any eigenforms $f \in S_{a+4}$, $g \in S_{b+c+4}$ and $h \in S_{b-c+2}$. If for a ``large'' prime $\ell$ in $\QQ_{f,g}$ and $s\geq 1$ the number $\ell^s$ divides the critical value $\Lambda'(\mathrm{Sym}^2(f) \otimes g \otimes h,a+b+5)$, then there is an eigenform $F \in S_{a-b,b-c,c+4}$ such that
\begin{equation} \label{eq-secondg3cong2} 
\lambda_{q}(F) \equiv_{\ell^s} \lambda_{q}(f)\bigl(\lambda_{q}(g)+q^{c+1} \lambda_{q}(h)\bigr)
\end{equation}
holds in $\QQ_{f,g,h,F}$ for all prime powers $q$.
\end{conjecture}

In this case, Mellit computed the critical value $\Lambda'(\Sym^2(f) \otimes g \otimes h,a+b+5)$ when $(a,b,c)=(16,14,0)$, for which conjecturally $s^{\gen}_{a-b,b-c,c+4}=1$. It is divisible by $\ell=71$ and indeed $N_q(a,b,c)$ is congruent, modulo $\ell$, to the right hand side of \eqref{eq-secondg3cong2} for all $q \leq 17$.

In the regular case, the endoscopic contribution to $e_c(\A_3,\VV_{a,b,c})$ consists conjecturally of two pieces, see Conjecture~\ref{conj-g3-reg}. Hence, in our picture, there is possibly a congruence connected to the other piece. But we have no numerical evidence, because $s^{\gen}_{a-b,b-c,c+4}$ is conjecturally larger than $1$ in all cases when the congruence could appear.

\begin{table}
\vbox{
\centerline{\def\quad{\hskip 0.3em\relax}
\vbox{\offinterlineskip
\hrule
\halign{&\vrule#& \quad \hfil#\hfil \strut \quad \cr
height2pt&\omit&&\omit&&\omit&&\omit&\cr
& $(a,b,c)$ && $e_c({\mathcal A}_3,\VV_{a,b,c})$ && $(a,b,c)$ && $e_c({\mathcal A}_3,\VV_{a,b,c})$ &\cr
\noalign{\hrule}
height2pt&\omit&&\omit&&\omit&&\omit&\cr
  & $(0,0,0)$ &&  $L^6+L^5+L^4+L^3+1$ && $(2,0,0)$ &&  $-L^3-L^2$ & \cr
  & $(1,1,0)$ &&  $-L$ && $(4,0,0)$ &&  $-L^3-L^2$ & \cr
  & $(3,1,0)$ &&  $0$ && $(2,2,0)$ &&  $0$ & \cr
  & $(2,1,1)$ &&  $1$ && $(6,0,0)$ &&  $-2L^3-L^2 $ & \cr
  & $(5,1,0)$ &&  $-L^4 $ && $(4,2,0)$ &&  $-L^5+L $ & \cr
  & $(4,1,1)$ &&  $1 $ && $(3,3,0)$ &&  $L^7-L $ & \cr
  & $(3,2,1)$ &&  $0 $ && $(2,2,2)$ &&  $1 $ & \cr
  & $(8,0,0)$ &&  $-L^3-L^2+S[12] $ && $(7,1,0)$ &&  $-L $ & \cr
  & $(6,2,0)$ &&  $L $ && $(6,1,1)$ &&  $-L^2+1 $ & \cr
  & $(5,3,0)$ &&  $0 $ && $(5,2,1)$ &&  $0 $ & \cr
  & $(4,4,0)$ &&  $0 $ && $(4,3,1)$ &&  $0 $ & \cr
  & $(4,2,2)$ &&  $L^4 $ && $(3,3,2)$ &&  $-L^6+1 $ & \cr
  & $(10,0,0)$ &&  $-2L^3-S[12]L^3 $ && $(9,1,0)$ &&  $-L^4+1 $ & \cr 
  & &&  $ -L^2+L $ && && & \cr
  & $(8,2,0)$ &&  $-L^5+L $ && $(8,1,1)$ &&  $1 $ & \cr
  & $(7,3,0)$ &&  $-L^6-L $ && $(7,2,1)$ &&  $0 $ & \cr
  & $(6,4,0)$ &&  $-L^7+L $ && $(6,3,1)$ &&  $-L^2 $ & \cr
  & $(6,2,2)$ &&  $ 0 $ && $(5,5,0)$ &&  $ L^9-L $ & \cr
  & $(5,4,1)$ &&  $0 $ && $(5,3,2)$ &&  $-L^3 $ & \cr
  & $(4,4,2)$ &&  $0 $ && $(4,3,3)$ &&  $-L^4+1 $ & \cr
  & $(12,0,0)$ &&  $-2L^3-L^2+S[16] $ && $(11,1,0)$ &&  $-L^4-S[12]L^4 $ & \cr
  & $(10,2,0)$ &&  $-L^5+L $ && $(10,1,1)$ &&  $-L^2-S[12]L^2+2 $ & \cr
  & $(9,3,0)$ &&  $-L^6+1 $ && $(9,2,1)$ &&  $0 $ & \cr
  & $(8,4,0)$ &&  $-L^7+L $ && $(8,3,1)$ &&  $-S[12] $ & \cr
  & $(8,2,2)$ &&  $L^4+S[12] $ && $(7,5,0)$ &&  $-L^8-L $ & \cr
  & $(7,4,1)$ &&  $0 $ && $(7,3,2)$ &&  $L^5 $ & \cr
  & $(6,6,0)$ &&  $L^{10}+S[0,10] $ && $(6,5,1)$ &&  $-L^2 $ & \cr
  & $(6,4,2)$ &&  $L^6-1 $ && $(6,3,3)$ &&  $1 $ & \cr
  & $(5,5,2)$ &&  $-L^8-L^3+1 $ && $(5,4,3)$ &&  $0 $ & \cr
  & $(4,4,4)$ &&  $-L^6+1 $ && $(14,0,0)$ &&  $-S[16]L^3-2L^3$ & \cr 
  &  &&  &&  &&  $-L^2+L+S[18] $ & \cr
  & $(13,1,0)$ &&  $-L^4-L-S[16]L+1 $ && $(12,2,0)$ &&  $-L^5-S[12]L^5+2L $ & \cr
  & $(12,1,1)$ &&  $-L^2+1 $ && $(11,3,0)$ &&  $-L^6-L $ & \cr
  & $(11,2,1)$ &&  $0 $ && $(10,4,0)$ &&  $-L^7+L+S[6,8] $ & \cr
  & $(10,3,1)$ &&  $-L^2-S[12]L^2+1 $ && $(10,2,2)$ &&  $ L^4 $ & \cr
  & $(9,5,0)$ &&  $-L^8+1 $ && $(9,4,1)$ &&  $0 $ & \cr
  & $(9,3,2)$ &&  $L^5-L^3 $ && $(8,6,0)$ &&  $-L^9+L $ & \cr
  & $(8,5,1)$ &&  $-S[12] $ && $(8,4,2)$ &&  $L^6-1 $ & \cr
  & $(8,3,3)$ &&  $-L^4+1 $ && $(7,7,0)$ &&  $L^{11}-2L $ & \cr
  & $(7,6,1)$ &&  $0 $ && $(7,5,2)$ &&  $L^7 $ & \cr
  & $(7,4,3)$ &&  $ 0 $ && $(6,6,2)$ &&  $-L^9+L^3+S[0,10] $ & \cr
  & $(6,5,3)$ &&  $ L^4 $ && $(6,4,4)$ &&  $0 $ & \cr
  & $(5,5,4)$ &&  $-L^8+1 $ && $ $ &&  $ $ & \cr
} \hrule}
}}
\end{table}

\begin{table}
\vbox{
\centerline{\def\quad{\hskip 0.3em\relax}
\vbox{\offinterlineskip
\hrule
\halign{&\vrule#& \quad \hfil#\hfil \strut \quad \cr
height2pt&\omit&&\omit&&\omit&&\omit&\cr
& $(a,b,c)$ && $e_c({\mathcal A}_3,\VV_{a,b,c})$ && $(a,b,c)$ && $e_c({\mathcal A}_3,\VV_{a,b,c})$ &\cr
\noalign{\hrule}
height2pt&\omit&&\omit&&\omit&&\omit&\cr
  & $(16,0,0)$ &&  $-S[18]L^3-2L^3-L^2$ && $(15,1,0)$ &&  $-L^4-S[16]L^4 $ & \cr
  & && $+L+S[20] $  && && $-S[18]L+1 $ & \cr  
  & $(14,2,0)$ &&  $-L^5+L+S[12,6] $ && $(14,1,1)$ &&  $-L^2-S[16]L^2+2 $ & \cr
  & $(13,3,0)$ &&  $-L^6-S[12]L^6-S[16]L+1 $ && $(13,2,1)$ &&  $0 $ & \cr
  & $(12,4,0)$ &&  $-L^7+L+S[8,8] $ && $(12,3,1)$ &&  $-L^2-S[16] $ & \cr
  & $(12,2,2)$ &&  $L^4+S[12]L^4-1+S[16] $ && $(11,5,0)$ &&  $-L^8-L $ & \cr
  & $(11,4,1)$ &&  $0 $ && $(11,3,2)$ &&  $L^5-L^3-S[12]L^3+1 $ & \cr
  & $(10,6,0)$ &&  $-L^9+L+S[4,10] $ && $(10,5,1)$ &&  $-L^2-S[12]L^2+1 $ & \cr
  & $(10,4,2)$ &&  $L^6+S[6,8]-1 $ && $(10,3,3)$ &&  $-L^4+1 $ & \cr
  & $(9,7,0)$ &&  $-L^{10}+1 $ && $(9,6,1)$ &&  $0 $ & \cr
  & $(9,5,2)$ &&  $ L^7-L^3 $ && $(9,4,3)$ &&  $0 $ & \cr
  & $(8,8,0)$ &&  $L^{12}+L+S[12]L+S[0,12] $ && $(8,7,1)$ &&  $L^2+S[12]L^2-S[12] $ & \cr
  & $(8,6,2)$ &&  $L^8+L^3+S[12]L^3-1 $ && $(8,5,3)$ &&  $S[12]L^4-S[12] $ & \cr
  & $(8,4,4)$ &&  $-S[12]L^6+S[12]+S[4,0,8] $ && $(7,7,2)$ &&  $-L^{10}+1 $ & \cr
  & $(7,6,3)$ &&  $0 $ && $(7,5,4)$ &&  $L^7-L^5 $ & \cr
  & $(6,6,4)$ &&  $-L^9+S[0,10] $ && $(6,5,5)$ &&  $-L^6+1 $ & \cr
  & $(18,0,0)$ &&  $-2S[20]L^3-3L^3-L^2$ && $(17,1,0)$ &&  $-2L^4-2S[18]L^4$ & \cr 
   & && $+L+S[22]$  && && $-S[20]L+1 $ & \cr
 & $(16,2,0)$ &&  $-2L^5-2S[16]L^5+2L $ && $(16,1,1)$ &&  $ -L^2-S[18]L^2+2 $ & \cr
  & $(15,3,0)$ &&  $-2L^6-L-S[18]L$ && $(15,2,1)$ &&  $0 $ & \cr
  & && $+1+S[12,7] $  && &&  & \cr 
  & $(14,4,0)$ &&  $-2L^7-2S[12]L^7+2L $ && $(14,3,1)$ &&  $-L^2-S[16]L^2$ & \cr
 & && && && $+1-S[18] $ & \cr
  & $(14,2,2)$ &&  $L^4+S[18]+S[12,6] $ && $(13,5,0)$ &&  $-2L^8-L-S[16]L $ & \cr
  & && && && $+1+S[8,9]$ & \cr
  & $(13,4,1)$ &&  $0 $ && $(13,3,2)$ &&  $L^5+S[12]L^5-L^3-1 $ & \cr
  & $(12,6,0)$ &&  $-2L^9+L+S[6,10] $ && $(12,5,1)$ &&  $-L^2-S[16] $ & \cr
  & $(12,4,0)$ &&  $ L^6+S[8,8]-1$ && $(12,3,3)$ &&  $-L^4-S[12]L^4+2 $ & \cr
  & $(11,7,0)$ &&  $-2L^{10}-L $ && $(11,6,1)$ &&  $0 $ & \cr
  & $(11,5,2)$ &&  $L^7-L^3-S[12]L^3$ && $(11,4,3)$ &&  $0 $ & \cr
  & && $+1 +S[6,3,6]$  && &&  & \cr  
  & $(10,8,0)$ &&  $-2L^{11}+2L $ && $(10,7,1)$ &&  $-S[12]L^2+1 $ & \cr
  & $(10,6,2)$ &&  $L^8+L^3+S[4,10]-1 $ && $(10,5,3)$ &&  $0 $ & \cr
  & $(10,4,4)$ &&  $S[6,8] $ && $(9,9,0)$ &&  $2L^{13}-L+S[12]+1 $ & \cr
  & $(9,8,1)$ &&  $0 $ && $(9,7,2)$ &&  $L^9-L^3 $ & \cr
  & $(9,6,3)$ &&  $S[3,3,7] $ && $(9,5,4)$ &&  $L^7-L^5 $ & \cr
  & $(8,8,2)$ &&  $-L^{11}+S[0,12] $ && $(8,7,3)$ &&  $-L^4-S[12] $ & \cr
  & $(8,6,4)$ &&  $ L^8-1 $ && $(8,5,5)$ &&  $-L^6+1 $ & \cr
  & $(7,7,4)$ &&  $-L^{10}-L^5+1 $ && $(7,6,5)$ &&  $0 $ & \cr
  & $(6,6,6)$ &&  $-L^9-L^8+1+S[0,10] $ && $ $ &&  $ $ & \cr
} \hrule}
}}
\end{table}


\begin{thebibliography}{99} 


\bibitem{Andrianov} A.N. Andrianov, V.G. Zhuravl\"ev:
{\sl Modular forms and Hecke operators.}
Translated from the 1990 Russian original by Neal Koblitz. Translations of Mathematical Monographs, 145. American Mathematical Society, Providence, RI, 1995.

\bibitem{Arthur} J.\ Arthur: 
{\sl An introduction to the trace formula.} In: Harmonic analysis, the trace formula, and Shimura varieties, Clay Math. Proc., Amer. Math. Soc., Providence, RI, 2005, 1--263.
In: Harmonic analysis, the trace formula, and Shimura varieties, 1--263, Clay Math. Proc., 4, Amer. Math. Soc., Providence, RI, 2005. 

\bibitem{B1} J.\ Bergstr\"om: 
{\sl Cohomology of moduli spaces of curves of genus three via point counts.}
J. Reine Angew. Math. \textbf{622} (2008), 155--187.

\bibitem{B2} J.\ Bergstr\"om: 
{\sl Equivariant counts of points of the moduli spaces of pointed hyperelliptic curves.}
Doc. Math. \textbf{14} (2009), 259--296.

\bibitem{BFG} J.\ Bergstr\"om, C.\ Faber, and G.\ van der Geer:
{\sl Siegel modular forms of genus~$2$ and level~$2$: conjectures and cohomological computations.}
IMRN, Art. ID rnn 100 (2008).

\bibitem{BFG2} J. Bergstr\"om, C. Faber, and G. van der Geer:
{\em Teichm\"uller modular forms and the cohomology of local systems on $M_3$.} In preparation. 

\bibitem{JBvdG} J.\ Bergstr\"om, G.\ van der Geer:
{\sl The Euler characteristic of local systems on the moduli of curves and abelian varieties of genus three.}
J. Topol. \textbf{1} (2008), 651--662.

\bibitem{BvdG} G.\ Bini, G.\ van der Geer: 
{\sl The Euler characteristic of local systems on the moduli of genus~$3$ hyperelliptic curves.}
 Math.\ Ann. \textbf{332} (2005), 367--379.

\bibitem{Borel} A. \ Borel:
{\sl Automorphic $L$-functions.} 
In: Automorphic forms, representations and $L$-functions
(Corvallis, OR, 1977), 27--61,
Proc. Sympos. Pure Math., 33, Part 2, Amer. Math. Soc., Providence, RI, 1979.

\bibitem{C-R} G. Chenevier, D. Renard:
{\sl Level one algebraic cusp forms of classical groups of small ranks.}
{\tt arXiv:1207.0724}.

\bibitem{C-F} C.\ Consani, C.\ Faber:
{\sl On the cusp form motives in genus~$1$ and level~$1$.  }
In: 
Moduli spaces and arithmetic geometry,  297--314, Adv.\ Stud.\ Pure Math.
\textbf{45}, Math.\ Soc.\ Japan, Tokyo, 2006.

\bibitem{Fermat}
H.\ Darmon, F.\ Diamond, and R.\ Taylor:
{\sl Fermat's last theorem.}
In: {Current developments in mathematics}, 1--154,
Int. Press, Cambridge, MA, 1994.

\bibitem{Deligne} P.\ Deligne: {\sl Formes modulaires et repr\'esentations
$\ell$-adiques}, S\'eminaire Bourbaki 1968-1969, exp. \textbf{355}, 139--172,
Lecture Notes in Math. \textbf{179}, 
Springer-Verlag, Berlin, 1971.

\bibitem{SGA} P. \ Deligne: {\sl Rapport sur la formule de trace.} In: Cohomologie \'etale, SGA $4\tfrac12$, Lecture Notes in Math. \textbf{569}. Springer-Verlag, Berlin-New York, 1977.

\bibitem{Deligne-L} P. \ Deligne: {\sl Valeurs de fonctions $L$ et p\'eriodes d'int\'egrales.}
In: Automorphic forms, representations and $L$-functions
(Corvallis, OR, 1977), 313--346,
Proc. Sympos. Pure Math., 33, Part 2, Amer. Math. Soc., Providence, RI, 1979.

\bibitem{DeligneWII}
P.\ Deligne: 
{\sl La conjecture de {W}eil. {II}.}
Inst. Hautes \'Etudes Sci. Publ. Math. \textbf{52} (1980), 137--252.

\bibitem{Deligne-letters} P.\ Deligne: {\rm Private communication}, letters dated June 4 (2009) and June 9 (2009).

\bibitem{Dummigan-sym} N. \ Dummigan: {\it Symmetric square $L$-functions and Shafarevich-Tate groups.} Experiment. Math. \textbf{10} (2001), no. 3, 383--400. 

\bibitem{Dummigan-sym-II} N. \ Dummigan: {\it Symmetric square $L$-functions and Shafarevich-Tate groups. II.} Int. J. Number Theory \textbf{5} (2009), no. 7, 1321--1345. 

\bibitem{DIK} N. Dummigan, T. Ibukiyama, and H. Katsurada: {\it
Some Siegel modular standard $L$-values, and Shafarevich-Tate groups}.
J. Number Theory \textbf{131} (2011), 1296--1330.

\bibitem{FvdG} C.\ Faber, G.\ van der Geer:
{\sl Sur la cohomologie des syst\`{e}mes locaux sur les espaces des modules des courbes de genre $2$ et des surfaces ab\'{e}liennes},
I, II. C.R.\ Acad.\ Sci.\ Paris, S\'er.\ I \textbf{338} (2004), 381--384, 467--470.

\bibitem{Faltings} G.\ Faltings:
{\sl On the cohomology of locally symmetric Hermitian spaces}.
Lecture Notes in Math. \textbf{1029}, Springer-Verlag, Berlin, 1983.

\bibitem{F-C} G.\ Faltings, C-L.\ Chai:
{\sl Degeneration of abelian varieties.}
Ergebnisse der Mathematik und ihrer Grenzgebiete (3), 22. Springer-Verlag, Berlin, 1990.

\bibitem{Freitag} E. \ Freitag:
{\sl Siegelsche Modulfunktionen.}
Grundlehren der Mathematischen Wissen\-schaften, 254. Springer-Verlag, Berlin, 1983.

\bibitem{FH} W.\ Fulton, J.\ Harris: {\sl Representation Theory.
A First Course}. Graduate Texts in Mathematics, 129.
Springer-Verlag, New York, 1991.

\bibitem{vdG1} G.\ van der Geer: {\sl
Siegel modular forms and their applications.}
In: J.\ Bruinier, G.\ van der Geer, G.\ Harder, D.\ Zagier:
The 1-2-3 of modular forms. Springer-Verlag, 2008.

\bibitem{vdG2} G.\ van der Geer:
{\sl Rank one Eisenstein cohomology of local systems on the moduli space of abelian varieties.} Sci. China Math. \textbf{54} (2011), no. 8, 1621--1634. 

\bibitem{Getz} E.\ Getzler: 
{\sl Topological recursion relations in genus~{$2$}.}
In: Integrable systems and algebraic geometry (Kobe/Kyoto, 1997), 
73--106. World Sci. Publishing, River Edge, NJ, 1998.

\bibitem{Gross-Savin} B.H. Gross, G. Savin: 
{\sl Motives with Galois group of type $G_2$: an exceptional the\-ta-correspondence.}
Compositio Math. \textbf{114} (1998), no. 2, 153--217.

\bibitem{Grundh}
C. Grundh: {\em
Computations of vector valued Siegel modular forms}.
Ph.D.~thesis, KTH, Stockholm, in preparation.

\bibitem{Haberland} K. Haberland: {\it Perioden von Modulformen einer Variablen und Gruppencohomologie.}
Math. Nachr. \textbf{112} (1983), 245--315.

\bibitem{Hain} R. \ Hain: {\it The rational cohomology ring of the moduli space of abelian 3-folds.}
Math. Res. Lett. \textbf{9} (2002), no. 4, 473--491. 

\bibitem{Ha1} G. Harder:
{\sl Eisensteinkohomologie und die Konstruktion gemischter Motive.}
Lecture Notes in Math. {\textbf 1562}. Springer-Verlag, Berlin, 1993.

\bibitem{Ha} G.\ Harder: {\sl A congruence between a Siegel and an
elliptic modular form.} In:
J.\ Bruinier, G.\ van der Geer, G.\ Harder, D.\ Zagier:
The 1-2-3 of modular forms. Springer-Verlag, 2008.

\bibitem{Ha2} G. Harder: {\it The Eisenstein motive for the cohomology of 
$\GSp_2(\ZZ)$.} 
In: Geometry and Arithmetic, 143--164,
EMS Series of Congress Reports, European Mathematical Society, Z\"urich, 2012.

\bibitem{Igusa} J. \ Igusa: {\sl On Siegel modular forms genus two. II.} 
Amer. \ J. \ Math. \textbf{86} (1964), 392--412.

\bibitem{Ikeda} T.\ Ikeda: {\sl Pullback of the lifting of elliptic cusp forms
and Miyawaki's conjecture.} Duke Math.J. \textbf{131} (2006), 469--497.

\bibitem{Katsurada-letter} H.\ Katsurada: {\rm Private communication}, letter dated August 19 (2011).

\bibitem{Katsurada-SK} H.\ Katsurada: {\sl Congruence of Siegel modular forms and special values of their standard zeta functions.} Math. Z. \textbf{259} (2008), no. 1, 97--111. 

\bibitem{Katsurada-Mizumoto} H.\ Katsurada, S. Mizumoto: {\sl Congruences of Hecke eigenvalues of Siegel modular forms.} Preprint (2011). 

\bibitem{Kottwitz} R.E. \ Kottwitz:
{\sl Points on some {S}himura varieties over finite fields.} J. Amer. Math. Soc. \textbf{5} (1992), 373--444.

\bibitem{Laumon} G. \ Laumon: 
{\sl Fonctions z\^etas des vari\'et\'es de {S}iegel de dimension trois.} Ast\'erisque \textbf{302} (2005), 1--66.

\bibitem{Serre-appendix} K. \ Lauter:
{\sl Geometric methods for improving the upper bounds on the number of rational points on algebraic curves over finite fields.} With an appendix by J.-P. Serre.
J. Algebraic Geom. \textbf{10} (2001), no. 1, 19--36.

\bibitem{Miyawaki} I.\ Miyawaki: {\sl Numerical examples of Siegel cusp
forms of degree $3$ and their zeta functions.} Memoirs of the
Faculty of Science, Kyushu University, Ser.\ A {\bf 46} (1992), 307--339.
 
\bibitem{Mizumoto} S. Mizumoto: 
{\sl Congruences for eigenvalues of Hecke operators on Siegel modular forms of degree two.}
Math. Ann. \textbf{275} (1986), no. 1, 149--161.

\bibitem{Morel} S. \ Morel:
{\sl On the cohomology of certain noncompact {S}himura varieties.} 
Annals of Mathematics Studies \textbf{173}, Princeton University Press, Princeton, NJ, 2010. 


\bibitem{Ram-Sha} D. \ Ramakrishnan, F. \ Shahidi: {\it Siegel modular forms of genus~2 attached to elliptic curves.} 
Math. Res. Lett. \textbf{14} (2007), no. 2, 315--332.

\bibitem{Satoh} T. \ Satoh: {\it On certain vector valued Siegel modular forms of degree two.}
Math. Ann. \textbf{274} (1986), no. 2, 335--352.

\bibitem{Scholl} A.J.\ Scholl: {\sl
Motives for modular forms.}  Invent.\ Math.
\textbf{100}  (1990),  no. 2, 419--430.

\bibitem{Serre-galmot} J.-P. \ Serre:
{\sl Propri\'et\'es conjecturales des groupes de Galois motiviques et des repr\'esentations $\ell$-adiques.} 
In:
Motives (Seattle, WA, 1991), 377--400,
Proc. Sympos. Pure Math., 55, Part 1, Amer. Math. Soc., Providence, RI, 1994.

\bibitem{Skoruppa} N.P. \ Skoruppa: {\it Computations of Siegel modular forms of genus two.} 
Math. Comp. \textbf{58} (1992), no. 197, 381--398. 

\bibitem{Stein} W.\ Stein: 
%  http://modular.fas.harvard.edu/tables/
The Modular Forms Database. \hfill\par
{\tt http://modular.math.washington.edu/Tables/}.

\bibitem{Taylor-Siegel3} R. \ Taylor: {\it On the $\ell$-adic cohomology of Siegel threefolds.}
Invent. Math. \textbf{114} (1993), no. 2, 289--310. 

\bibitem{Tsushima} R.\ Tsushima:
{\sl An explicit dimension formula for the spaces of generalized
automorphic forms with respect to ${\rm Sp}(2,\ZZ)$.}
Proc.\ Japan Acad.\ Ser.\ A Math.\ Sci. \textbf{59} (1983), no. 4, 139--142.

\bibitem{Tsuyumine} S.\ Tsuyumine: 
{\sl On Siegel Modular Forms of Degree Three}.
Am. \ J. \ Math. \textbf{108} (1986), no. 4, 755--862.

\bibitem{Weissauer-book} R. Weissauer: {\sl Endoscopy for GSp(4) and the cohomology of Siegel modular threefolds.} Lecture Notes in Math. \textbf{1968}. Springer-Verlag, Berlin, 2009.

\bibitem{Weissauer} R. Weissauer: 
{\sl The trace of Hecke operators on the space of classical holomorphic Siegel modular forms of genus two.}
{\tt arXiv:0909.1744}.

\bibitem{Zagier-mod} D. \ Zagier: {\it Modular forms whose Fourier coefficients involve zeta-functions of quadratic fields.} Modular functions of one variable VI, 105--169,
Lecture Notes in Math. \textbf{627}. Springer-Verlag, Berlin, 1977.

\bibitem{Z} D.\ Zagier:
{\sl Sur la conjecture de Saito-Kurokawa (d'apr\`es H.\
Maass).}  Seminar on Number Theory,  Paris 1979--80,  371--394,
Progr.\ Math., 12, Birkh\"auser, Boston, MA, 1981.
\end{thebibliography}
\end{document}